\pdfoutput=1
\documentclass[toc]{mathprint} %
\usepackage{rutarmacros}
\usepackage{macros}

\title[Lower box dimensions]{Lower box dimension of infinitely generated self-conformal sets}

\author[Banaji]
  {Amlan Banaji}
  {Mathematical Sciences, Loughborough University, Loughborough, LE11 3TU, United Kingdom}
  {A.F.Banaji@lboro.ac.uk}

\author[Rutar]
  {Alex Rutar}
  {University of St Andrews, Mathematical Institute, St Andrews, KY16 9SS, Scotland}
  {alex@rutar.org}

\begin{document}
\begin{abstract}
    Let $\Lambda$ be the limit set of an infinite conformal iterated function system and let $F$ denote the set of fixed points of the maps.
    We prove that the box dimension of $\Lambda$ exists if and only if
    \begin{equation*}
        \overline{\dim}_{\mathrm B} F\leq \max \{\dim_{\mathrm H} \Lambda, \underline{\dim}_{\mathrm B} F\}.
    \end{equation*}
    In particular, this provides the first examples of sets of continued fraction expansions with restricted digits for which the box dimension does not exist.

    More generally, we establish an explicit asymptotic formula for the covering numbers $N_r(\Lambda)$ in terms of $\dim_{\mathrm{H}}\Lambda$ and the covering function $r\mapsto N_r(F)$, where $N_r(\cdot)$ denotes the least number of open balls of radius $r$ required to cover a given set.
    Such finer scaling information is necessary: in general, the lower box dimension $\underline{\dim}_{\mathrm B} \Lambda$ is not a function of the Hausdorff dimension of $\Lambda$ and the upper and lower box dimensions of $F$, and we prove sharp bounds for $\underline{\dim}_{\mathrm B} \Lambda$ in terms of these three quantities.
    \vspace{0.5cm}


    \vspace{0.2cm}

\end{abstract}

\section{Introduction}
\subsection{Dynamical invariance and box dimension}
One of the fundamental aspects of the study of dynamical systems is the structure of invariant sets.
More precisely, suppose $X$ is a set and $f\colon X\to X$ is a function.
We say that a set $\Lambda\subset X$ is \emph{invariant} if $f(\Lambda)=\Lambda$.
Invariant sets for dynamical systems, particularly when the dynamics are expanding or chaotic, often have highly intricate local structure.
As a result, a standard lens through which to view the geometry of invariant sets is that given by \emph{dimension theory}.
The dimension theory of dynamical systems has its origins in the seminal work of Bowen \cite{zbl:0439.30032} on quasicircles and Ruelle \cite{zbl:0506.58024} on conformal repellers.
It has since developed into a substantial field in its own right.
For an overview, see, for instance, the monographs \cite{zbl:0895.58033,zbl:1161.37001}, surveys \cite{zbl:1230.37032,zbl:1001.37018}, and notable recent progress \cite{zbl:1387.37026,zbl:1414.28014,zbl:1427.37018,zbl:1531.37017}.

Of the notions of fractal dimension, one of the most familiar, and perhaps the simplest to define, is the \emph{box dimension}, also known as \emph{Minkowski dimension}.
More precisely, the lower and upper box dimensions of a non-empty bounded set $E \subset \R^d$ are defined respectively as
\begin{equation*}
    \dimlB E = \liminf_{r\to 0} \frac{\log N_r(E)}{\log (1/r)}, \qquad \dimuB E  = \limsup_{r\to 0} \frac{\log N_r(E)}{\log (1/r)}
\end{equation*}
where $N_r(E)$ denotes the least number of open balls of radius $r$ required to cover $E$.
If these two notions coincide, we call the common value the box dimension of $E$ and denote it by $\dimB E$.
If the box dimension does not exist, this indicates that the extent to which the set $E$ `fills up space' varies substantially at different scales.

It is well known that for general sets $E$, the inequalities $\dimH E \leq \dimlB E \leq \dimuB E$ always hold and moreover can be strict, where $\dimH E$ denotes the Hausdorff dimension of $E$.
On the other hand, if $E$ is invariant for a map $f$ which is uniformly expanding on $E$, then in many cases these will be equalities.
A core question in the dimension theory of dynamical systems, and the question motivating this paper, is the following.
\begin{question}
    For which sets that are invariant for an expanding dynamical system does the box dimension exist?
\end{question}
It was proved independently in \cite{zbl:0876.58039,zbl:0862.58042} (generalising a result in \cite{zbl:0683.58034}) that if $E$ is a compact subset of a Riemannian manifold which is invariant under a smooth conformal map $f$ that is uniformly expanding on $E$, then $\dimH E = \dimlB E = \dimuB E$.
Moreover, if $E$ is a self-similar or self-conformal set (these are the most familiar classes of fractal sets), then $\dimH E = \dimlB E = \dimuB E$.

The assumption that the dynamics are conformal (meaning that the differential $f'_x$ is a scalar multiple of a similarity map at each $x \in E$), and the assumption that the invariant set $E$ is compact, are important.
It has been known that a gap between Hausdorff and lower box dimension is possible if $f$ is non-conformal since the work of Bedford \cite{Bedford} and McMullen \cite{zbl:0539.28003} in 1984 (the invariant sets are self-affine Bedford--McMullen carpets).
The box and Hausdorff dimension can also differ for Julia sets of certain non-rational hyperbolic functions \cite{zbl:1034.30019,zbl:1050.30017}.
Moreover, from work of Mauldin \& Urbański \cite{zbl:0852.28005,zbl:0940.28009} in the 1990s one can see that there are non-compact sets which are invariant for a dynamical system given by a function which extends to a smooth expanding conformal map, and whose Hausdorff and box dimensions differ.
In the latter case, the invariant set is an infinitely generated self-conformal set; such sets are described in detail below.

On the other hand, there has been much less progress on establishing a gap between lower and upper box dimension (i.e.\ non-existence of box dimension).
In the non-conformal setting, it is a major problem to determine whether the box dimension of every self-affine set exists, and while this is known to be true in many cases \cite{Bedford,zbl:0539.28003,zbl:0642.28005,zbl:1414.28014}, the general problem remains open.
The usual examples of sets whose box dimension may not exist, such as certain sequences, Cantor-like sets \cite[§2]{zbl:1285.28011}, and inhomogeneous attractors \cite{zbl:1342.28013}, are not dynamically invariant.
The only non-existence result for dynamically invariant sets of which we are aware is due to Jurga \cite{zbl:1531.37017}, who recently showed that there exists a compact subset of the torus which is invariant under a smooth expanding toral endomorphism and whose box dimension does not exist.
Crucial to Jurga's example is the two-dimensional nature of the torus and the non-conformal nature of the dynamics.

In this paper, we give the first examples (as far as we are aware) of \emph{non-compact} sets whose box dimension does not exist and which are invariant for a dynamical system given by a function which extends to a smooth uniformly expanding \emph{conformal} map.
This is the case even for a particularly old and well-studied conformal map: the \emph{Gauss map} $g\colon(0,1]\to S^1$ given by $g(x)=1/x$, where $S^1$ denotes the circle $\R\!/\!\Z$.
Canonical examples of invariant sets for the Gauss map are the numbers with continued fraction expansions having digits restricted to a given subset of $\N$.
This is also the case for the set $\Lambda$ in the theorem below.
\begin{itheorem}\label{it:gauss-invariant}
    There exists a Borel set $\Lambda\subset(0,1)$ which is invariant under the Gauss map and whose box dimension does not exist.
\end{itheorem}

\cref{it:gauss-invariant} is in fact a consequence of our study of \emph{infinitely generated} self-conformal sets in this paper.
Such sets were first introduced in \cite{zbl:0852.28005}, and as mentioned above, motivating examples of infinitely generated self-conformal sets are real or complex numbers whose continued fraction expansions are restricted to a given countable set \cite{zbl:0940.28009}.
Infinitely generated self-conformal sets are particularly well-studied in the literature; for a certainly incomplete selection, see for instance \cite{zbl:07662347,zbl:07844477,zbl:1437.37032,zbl:1319.11050,zbl:1233.11084,zbl:1482.37002} and \cite[Section~9.2]{zbl:1467.28001} and more references therein.
Infinite systems are also useful in the presence of non-uniformly expanding dynamics or parabolicity, since one can often associate with such a system an `induced' infinite but uniformly expanding system \cite{zbl:0982.37045,zbl:1013.28007}.
Furthermore, infinite systems have been used to calculate dimensions of sets which are important in complex dynamics \cite{zbl:1500.37035}.
While the Hausdorff and upper box dimensions of infinitely generated self-conformal sets are well-understood \cite{zbl:0852.28005,zbl:0940.28009}, despite over two decades since the original results, much less is known concerning the lower box dimension.

In this paper we remedy this gap by providing a precise formula for the lower box dimension of such sets; this formula is substantially more complicated than the one for upper box dimension.
Moreover, we characterise when the lower and upper box dimensions of such sets coincide.

\subsection{Countably generated self-conformal sets}\label{s:conformal}
In order to state our main result precisely, we introduce our main objects under consideration: countable conformal iterated function systems.

Following \cite{zbl:0852.28005}, let $X$ be a compact connected subset of $\R^d$ with the Euclidean norm and let $\mathcal{I}$ be a countable index set.
Fix a family of injective uniformly contracting maps $S_i \colon X \to X$ for $i \in \mathcal{I}$: that is, there is a constant $0<c<1$ so that for all $x,y\in X$ and $i\in\mathcal{I}$,
\begin{equation*}
    0 < \norm{S_i(x)-S_i(y)}\leq c\cdot \norm{S_i(x)-S_i(y)}.
\end{equation*}
We use symbolic notation on the set $\mathcal{I}^*$ of finite sequences on $\mathcal{I}$, equipped with the operation of concatenation.
Given $\gamma = (i_1,i_2,\dotsc) \in \mathcal{I}^{\N}$ and $n \in \N$, we write $\gamma|_{n} \coloneqq (i_1,\dotsc,i_n) \in \mathcal{I}^n$.
For $n\in\N\cup\{0\}$ and $\mtt{i}=(i_1,\ldots,i_n)\in\mathcal{I}^n$, we write
\begin{equation*}
    S_{\mtt{i}}=S_{i_1}\circ\cdots\circ S_{i_n}.
\end{equation*}
Then, associated with the IFS $\{S_i\}_{i\in\mathcal{I}}$ is the \emph{limit set} (also called the \emph{attractor}):
\begin{equation*}
    \limitset \coloneqq \bigcup_{\gamma \in \mathcal{I}^\mathbb{N}} \bigcap_{n=1}^\infty S_{\gamma\npre{n}}(X).
\end{equation*}
Equivalently, $\limitset$ is the largest subset of $X$ (by inclusion) satisfying the invariance relationship
\begin{equation*}
    \limitset = \bigcup_{i \in \mathcal{I}} S_i(\limitset).
\end{equation*}
Note that $\limitset$ is not in general a compact set.
On the other hand, if for each $x\in X$ there are only finitely many indices $i$ such that $x\in S_i(X)$ (which is the case for all systems which we will consider in this paper), then $\limitset$ is an $F_{\sigma\delta}$ subset of $X$.

Throughout, $\limitset$ will denote the limit set and $\fixedpts$ will denote the (countable) set of fixed points of the contractions $\{S_i\}_{i\in\mathcal{I}}$.
There is nothing special about this choice; in general, one could choose any point $x_i\in S_i(X)$ for $i\in\mathcal{I}$ and let $\fixedpts=\{x_i:i\in\mathcal{I}\}$.
More discussion can be found in \cref{ss:form}.

\begin{definition}\label{d:cifs}
    We say that the IFS $\{S_i\}_{i\in\mathcal{I}}$ is \emph{conformal} if the following additional properties are satisfied:
    \begin{enumerate}[r]
        \item\label{i:conformal} \emph{Conformality}:
            There exists an open, bounded, connected subset $V \subset \R^d$ such that $X \subset V$ and such that for each $i \in \mathcal{I}$, $S_i$ extends to a conformal $C^{1+\varepsilon}$ diffeomorphism on $V$.
        \item\label{i:bdp} \emph{Bounded distortion}:
            There exists $K\geq 1$ such that $\norm{S_{\mtt{i}}'(x)} \leq K\snorm{S_{\mtt{i}}'(y)}$ for all $x,y \in V$ and $\mtt{i} \in \mathcal{I}^*$.
            Here, $S_{\mtt{i}}'(x)$ denotes the Jacobian of the map $S_{\mtt{i}}$ at $x$ and $\norm{\cdot}$ denotes the spectral matrix norm.
    \end{enumerate}
\end{definition}
In light of the uniformly contracting property, if we define
\begin{equation*}
    \rho(\mtt{i}) = \sup_{x\in X}\norm{S_{\mtt{i}}'(x)}\qquad\text{for}\qquad\mtt{i}\in\mathcal{I}^*,
\end{equation*}
then $\xi\coloneqq\sup_{i\in\mathcal{I}}\rho(i) <1$.
Moreover, by the chain rule, sub-multiplicativity of the matrix norm, and the bounded distortion property, for any $\mtt{i},\mtt{j}\in\mathcal{I}^*$,
\begin{equation}\label{e:rho-approx-mul}
    K^{-1}\rho(\mtt{i})\rho(\mtt{j})\leq \rho(\mtt{i}\mtt{j})\leq \rho(\mtt{i})\rho(\mtt{j}).
\end{equation}
We also require some standard separation conditions; see, for instance, \cite{zbl:0852.28005}.
\begin{enumerate}[r]
    \setcounter{enumi}{2}
    \item\label{i:osc}\emph{Open set condition}:
        The set $X$ has non-empty topological interior $U$, and $S_i(U) \subset U$ for all $i \in \mathcal{I}$ and $S_i(U) \cap S_j(U) = \varnothing$ for all $i,j \in \mathcal{I}$ with $i \neq j$.
    \item\label{i:cone} \emph{Cone condition}:
        $\displaystyle\inf_{x \in X} \inf_{r \in (0,1)} r^{-d}\mathcal{L}^d (B(x,r) \cap U) > 0$.
\end{enumerate}
See \cref{s:separation} for more discussion concerning separation conditions.

Throughout this paper, we will assume that $\{S_i\}_{i\in\mathcal{I}}$ is a conformal IFS satisfying the open set condition and the cone condition.
We will refer to such a system in shorthand as a \cifs.

In order to study the dimension theory of the limit set $\limitset$ of a \cifs{}, Mauldin \& Urbański~\cite{zbl:0852.28005} defined the \emph{topological pressure} $\pressure\colon (0,\infty) \to [-\infty,\infty]$ by
\begin{equation}\label{e:MUpressure}
    \pressure(t) \coloneqq \lim_{n \to \infty} \frac{1}{n} \log\sum_{\mtt{i}\in \mathcal{I}^n} \rho(\mtt{i})^t.
\end{equation}
The limit necessarily exists by a sub-additivity argument using \cref{e:rho-approx-mul}.
In the same paper, they established a formula for the Hausdorff dimension of the limit set in terms of the pressure:
\begin{equation}\label{e:hausdorff-formula}
    \dimH\limitset = \inf\{t \geq 0 : \pressure(t) < 0\}.
\end{equation}
We also know from \cite{zbl:0940.28009} that the upper box and packing dimensions coincide and are given by
\begin{equation*}
    \dimuB \limitset = \dimP \limitset = \max\{\dimH \limitset, \dimuB \fixedpts\},
\end{equation*}
see \cite{zbl:0940.28009,zbl:07662347}.
Other notions of dimension, such as the upper intermediate dimensions \cite{zbl:07662347} and Assouad-type dimensions \cite{zbl:07844477}, have also been studied.

\subsection{Formula for the lower box dimensions}
The lower box dimension of general sets has some interesting properties not shared by Hausdorff or upper box dimension.
For instance, the lower box dimension need not be finitely stable: it is not too challenging to construct sets $E_1,E_2$ such that $\dimlB (E_1\cup E_2)>\max\{\dimlB E_1,\dimlB E_2\}$ by choosing the scales at which each set is ``large'' to be very different.
As we will see, this sensitivity of the lower box dimension to the finer scaling properties of the underlying set leads to much more interesting behaviour for the limit set of a \cifs{} than appears for the upper box dimension.

For the limit set of a \cifs, the following bounds for lower box dimension are immediate from the work of Mauldin \& Urbański:
\begin{equation}\label{e:trivialbounds}
    \max \{\dimH \limitset, \dimlB \fixedpts\} \leq \dimlB \limitset \leq \dimuB\limitset = \max \{\dimH \limitset, \dimuB \fixedpts\}.
\end{equation}
In general, the lower bound in \cref{e:trivialbounds} is sharp in a sense which will become clear below.
In contrast, it turns out the upper bound is not sharp in general.
In fact, our main result provides a precise classification of the existence of the box dimension of $\limitset$ which equivalently states that the box dimension of $\limitset$ exists if and only if the bounds in \cref{e:trivialbounds} coincide.
\begin{itheorem}\label{it:box-exist}
    Let $\limitset$ be the limit set of a \cifs.
    Then $\dimlB\limitset=\dimuB\limitset$ if and only if
    \begin{equation*}
        \dimuB\fixedpts\leq \max \{\dimH \limitset, \dimlB \fixedpts\}.
    \end{equation*}
\end{itheorem}
In particular, non-existence of the box dimension is common.

\cref{it:box-exist} will follow from an explicit asymptotic formula for $N_r(\limitset)$ in terms of the scaling function $r\mapsto N_r(\fixedpts)$ and the Hausdorff dimension of the limit set $\limitset$, which we now state.

Given the set of fixed points $F$ and $r\in(0,1)$, we define the \emph{box dimension estimate at scale $r$} by
\begin{equation*}
    s(r)\coloneqq \frac{\log N_r(F)}{\log(1/r)}.
\end{equation*}
For $0<r<1$ and $0<\theta \leq 1$, we define
\begin{equation}\label{e:psi-def}
    \Psi(r,\theta)\coloneqq (1-\theta)\dimH \limitset + \theta s(r^\theta)
\end{equation}
and set $\Psi(r,0)\coloneqq \lim_{\theta\to 0}\Psi(r,\theta)$.
We then set
\begin{equation*}
    \psi(r)\coloneqq \max_{\theta\in[0,1]}\Psi(r,\theta);
\end{equation*}
the maximum exists by upper semi-continuity of the map $\theta\mapsto N_{r^\theta}(F)$.
A depiction of the function $\psi(r)$ in terms of the box dimension estimate $s(r)$ can be found in \cref{f:threshold}.
\begin{figure}[t]
    \centering
    \begin{tikzpicture}[>=stealth,xscale=3,yscale=9]
    \draw[->] (-0.1,0) -- (4.2,0) node[below left]{$x=\log\log(1/r)$};
    \draw[->] (0,-0.033333) -- (0,0.56666666);
    \draw[thick, dotted] (-0.1,0.26) -- node[below]{$\dimH\limitset$} (4,0.26);
    \draw[thick, dashed] plot file {figures/threshold/original.txt};
    \draw[thick] plot file {figures/threshold/modified.txt};
\end{tikzpicture}
    \caption{A plot of the functions $s(r)$ (dashed), the function $\psi$ (solid), and $\dimH\limitset$ (dotted).
        The domain has been transformed by the order-reversing map $x=\log\log(1/r)$---see \cref{p:covering-class} for more detail.
    }
    \label{f:threshold}
\end{figure}
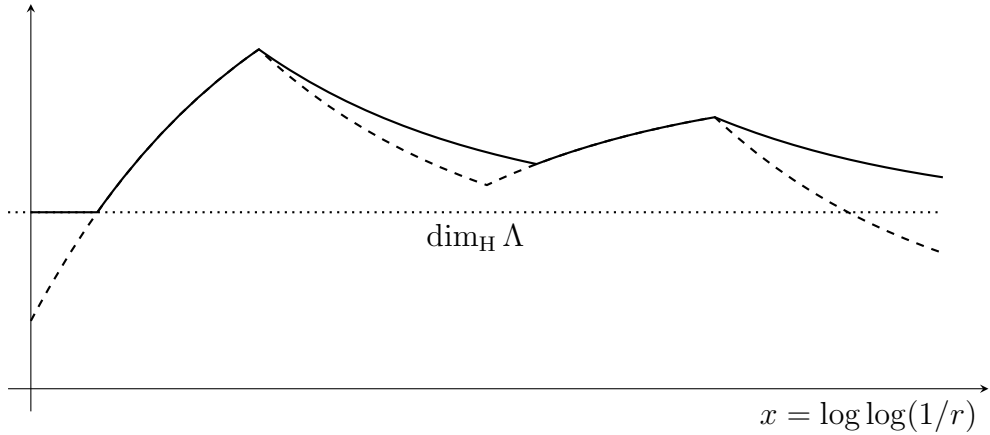
\begin{itheorem}\label{it:main}
    Let $\limitset$ be the limit set of a \cifs{} on $\R^d$ with fixed points $F$ and function $\psi$ as above.
    Then
    \begin{equation*}
        \lim_{r\to 0}\left(\frac{\log N_r(\limitset)}{\log(1/r)}-\psi(r)\right)=0.
    \end{equation*}
    In particular,
    \begin{equation*}
    	\dimlB\limitset=\liminf_{r\to 0}\psi(r).
    \end{equation*}
\end{itheorem}
We now make several comments on this result.
\begin{enumerate}[nl]
    \item The formula only depends on the contraction ratios through the Hausdorff dimension.
        On the other hand, it very much depends on the scaling properties of the set of fixed points.
    \item The formula can depend on the Hausdorff dimension even when $\dimH\limitset < \dimlB\fixedpts$.
    \item Heuristically, the formula says that $N_r(\limitset)\geq N_r(F)$ should be as small as possible while still being at least $\dimH\limitset$-dimensional between all pairs of scales.
        This heuristic can be made precise; the details can be found in \cref{ss:rate-intro}.
        In some sense, this is a more quantitative version of the observation that $\limitset$ contains Ahlfors--David $\lambda$-regular subsets for all $\lambda<\dimH\limitset$.
    \item Since $\lim_{r\to 0}\Psi(r,0)=\dimH\limitset$ and $\liminf_{r\to 0}\Psi(r,1)=\dimlB\fixedpts$, the trivial bound
        \begin{equation*}
            \dimlB\limitset=\liminf_{r\to 0}\psi(r)\geq\max\{\dimH\limitset,\dimlB F\}
        \end{equation*}
        corresponds to the endpoints of the estimate $\psi(r)$.
    \item A straightforward argument using only the definition of the upper box dimension gives that
        \begin{equation*}
            \dimuB\limitset=\limsup_{r\to 0}\psi(r) = \max\{\dimH\limitset, \dimuB\fixedpts\}.
        \end{equation*}
        See the proof of \cref{l:small-r-reg} for more detail.
        In particular, we recover the known results for the upper box dimension from \cite{zbl:0940.28009}.
\end{enumerate}

As an application of \cref{it:main}, we obtain the following bounds only assuming coarse scaling properties of $\fixedpts$.
Given $0\leq s\leq t\leq \alpha$ and $0\leq h\leq \alpha$, define the compact interval
\begin{equation*}
    \mathcal{D}(h,s,t,\alpha) =
    \begin{cases}
        \{h\} &: t\leq h,\\
        \left[\max\{h,s\}, h + \frac{(t-h) (\alpha - h) s}{\alpha\cdot t-h\cdot s}\right] &: t>h.
    \end{cases}
\end{equation*}
Observe that $\mathcal{D}(h,s,t,\alpha)$ is a not a singleton if and only if $0<h<t$ and $0<s<t$.
We then have the following result.
\begin{itheorem}\label{it:fixed}
    Let $\limitset$ be the limit set of a \cifs{} on $\R^d$.
    Then
    \begin{equation*}
        \dimlB\limitset\in\mathcal{D}\bigl(\dimH\limitset, \dimlB\fixedpts,\dimuB\fixedpts, d\bigr).
    \end{equation*}
    Moreover, these bounds are sharp in the following sense: for any $0<h<d$, $0\leq s\leq t\leq d$, and
    \begin{equation*}
        \beta\in\mathcal{D}(h,s,t,d),
    \end{equation*}
    there is a \cifs{} of similarity maps on $\R^d$ with limit set $\limitset$ and fixed points $\fixedpts$ such that $\dimH\limitset=h$, $\dimlB\fixedpts = s$, $\dimuB\fixedpts = t$, and $\dimlB\limitset=\beta$.
\end{itheorem}
In particular, we see that the inequalities $\dimH \limitset \leq \dimlB \limitset \leq \dimuB \limitset$ can be strict or non-strict in any combination.
This is in stark contrast to the case of finitely generated self-conformal sets, whose Hausdorff and upper box dimensions are always equal.
\begin{remark}\label{r:inf-bound}
    The family of sets $\mathcal{D}(h,s,t,\alpha)$ is monotonically increasing in $\alpha$ for fixed $h,s,t$.
    Moreover,
    \begin{align*}
        \mathcal{D}(h,s,t,\infty)&\coloneqq \lim_{\alpha\to\infty}\mathcal{D}(h,s,t,\alpha)\\*
                                 &=
        \begin{cases}
            \{h\} &: t\leq h,\\
            \left[\max\{h,s\},s + \left(1-\frac{s}{t}\right)h\right] &: t > h.
        \end{cases}
    \end{align*}
    This gives non-trivial bounds which are valid in all dimensions simultaneously.
    The right endpoint of $\mathcal{D}(h,s,t,\infty)$ is always at most $\max\{h,t\}$, and is equal to this upper bound if and only if the box dimension of $\limitset$ exists.
\end{remark}

\subsection{Reformulation in terms of growth rates of covering numbers}\label{ss:rate-intro}
To conclude, given a \cifs{} with limit set $\limitset$, we recast our asymptotic formula for $N_r(\limitset)$ in terms of a certain minimal growth rate of covering numbers.
The purpose of this is to make rigorous the claim following the statement of \cref{it:main} that $N_r(\limitset)\geq N_r(\fixedpts)$ is as small as possible while being at least $\dimH\limitset$-dimensional between all pairs of scales.

Let $d\in\N$ and fix a non-empty bounded subset $E\subset\R^d$.
Then for $0<r<1$, define
\begin{equation*}
    s_E(r)\coloneqq\frac{\log N_r(E)}{\log(1/r)}.
\end{equation*}

\begin{definition}
    Let $0\leq \lambda\leq\alpha$ be arbitrary.
    Let $\mathcal{G}(\lambda,\alpha)$ denote the set of continuous functions $g\colon\R\to[\lambda,\alpha]$ such that
    \begin{equation*}
        \diniu+ g(x)\in[\lambda-g(x),\alpha-g(x)],
    \end{equation*}
    where
    \begin{equation*}
    \diniu+ g(x)\coloneqq \limsup_{\varepsilon\to 0^+}\frac{g(x+\varepsilon)-g(x)}{\varepsilon}
    \end{equation*}
    is the upper right Dini derivative of $g$ at $x$.
\end{definition}
We say that two functions $f,g\colon\R\to\R$ are \emph{asymptotically equivalent}, and write $f\eqinf g$, if
\begin{equation*}
    \lim_{x\to \infty}\bigl(f(x)-g(x)\bigr)=0.
\end{equation*}
It is clear that $\eqinf$ is an equivalence relation.

The point here is that for a non-empty bounded set $E\subset\R^d$, after a change of domain, the function $s_E$ is asymptotically equivalent to a member of $\mathcal{G}(0,d)$.
This observation was essentially made in \cite {zbl:1509.28005}, and is made precise in the following proposition.
\begin{proposition}\label{p:covering-class}
    Let $E\subset\R^d$ be non-empty and bounded with associated function $f\colon\R\to\R$ defined by
    \begin{equation*}
        f(x)=s_E(\exp(-\exp(x)))
    \end{equation*}
    Then there exists $g\in\mathcal{G}(0,d)$ such that $f\eqinf g$.
\end{proposition}
The proof is straightforward though requires a slight amount of attention.
It follows from a slightly stronger version stated in \cref{p:covering-class-improved}.

With \cref{p:covering-class} in mind, we make the following definition.
\begin{definition}
    Let $E\subset\R^d$ be non-empty and bounded.
    We say that $E$ has \emph{covering class} $g\in\mathcal{G}(0,d)$ if $f\eqinf g$ where $f(x)=s_E(\exp(-\exp(x)))$.
\end{definition}
Of course, the covering class is only well-defined up to asymptotic equivalence.
In general, if $E$ has covering class $g$, then $\dimlB E=\liminf_{x\to\infty}g(x)$ and $\dimuB E=\limsup_{x\to\infty}g(x)$.
One can interpret the numbers $\lambda$ and $\alpha$ as saying that uniformly (in a weak exponential sense) over all pairs of scales and on average in space, $E$ is at least $\lambda$-dimensional and at most $\alpha$-dimensional.

By a straightforward argument given in \cref{p:relative-compact}, the class $\mathcal{G}(\lambda,\alpha)$ is closed under taking suprema and infima.
Recalling \cref{f:threshold}, we will use \cref{it:main} to prove the following result.
\begin{itheorem}\label{it:rate-form}
    Let $\limitset$ be the limit set of a \cifs{} on $\R^d$ with fixed points $F$.
    Let $F$ have covering class $f$ and let
    \begin{equation*}
        g\in\mathcal{G}\bigl(\dimH\limitset,d\bigr)
    \end{equation*}
    be the pointwise minimal function satisfying $f\leq g$.
    Then $\limitset$ has covering class $g$.
\end{itheorem}
\begin{remark}
    If $F$ has covering class $f\in\mathcal{G}(0,\alpha)$, then we in fact prove that $\limitset$ has covering class $g$ where $g\in\mathcal{G}\bigl(\dimH\limitset,\max\{\dimH\limitset, \alpha\}\bigr)$ is the pointwise minimal function satisfying $f\leq g$.
    Moreover, in this case, the proof of \cref{it:fixed} implies that the bound holds with $\alpha$ in place of $d$.
    In particular, by \cref{p:covering-class-improved},
    \begin{equation*}
        \dimlB\limitset\in\mathcal{D}(h,\dimlB\fixedpts,\dimuB\fixedpts,\dimqA\fixedpts),
    \end{equation*}
    where $\dimqA\fixedpts$ is the quasi-Assouad dimension of $\fixedpts$ introduced in \cite{zbl:1345.28019}.
\end{remark}

\subsection{Further research}
There are several possible directions for future research.
It is likely that methods in this paper would also extend to the case of the lower intermediate dimensions of the limit set of a \cifs{}, which would give a complete answer to \cite[Question~3.7]{zbl:07662347}.
In another direction, little is known about the box dimensions of sets generated by countably many \emph{affine} contractions, despite interest in the Hausdorff dimension of such sets \cite{arxiv:2405.00520,zbl:1292.28016} (some other work in this direction can be found in \cite{arxiv:2404.10749}).
Furthermore, it may be of interest to study existence of the  box dimension for the many possible generalisations of a CIFS, such as random, graph-directed, and/or overlapping systems \cite{zbl:1353.37111,zbl:1412.37032,zbl:1225.37063}.

We also note a connection with the study of \emph{inhomogeneous self-conformal sets}.
As established by Fraser \cite{zbl:1342.28013}, the lower box dimension of the attractor of an inhomogeneous self-similar IFS also depends on the covering properties of the condensation set at different scales through the \emph{covering regularity exponent}.
It seems plausible to the authors that the techniques in this paper would be useful in the study of box dimensions of inhomogeneous self-conformal sets.

\subsection{Notation, conventions, and structure of paper}
Throughout the paper, we work in $\R^d$ for $d\in\N$ equipped with the Euclidean norm.
The ball $B(x,r)$ is the open ball centred at $x$ with radius $r$.
The covering number $N_r(E)$ is the least number of open balls of radius $r$ required to cover $E$.

We will also find it useful to use asymptotic notation.
Given a set $A$ and functions $f,g\colon A\to\R$, we write $f\lesssim g$ if there is a constant $C>0$ such that $f(a)\leq C g(a)$ for all $a\in A$.
We write $f\approx g$ if $f\lesssim g$ and $f\gtrsim g$.
The constants in the asymptotic notation will always be allowed to implicitly depend on the underlying IFS; any other dependence will be explicitly indicated by a subscript, such as $\lesssim_\varepsilon$.

In \cref{s:proof} we prove our key result \cref{it:main} giving the asymptotic formula for covering numbers.
In \cref{s:consequences} we prove several consequences of \cref{it:main}, including the bounds in the first half of \cref{it:fixed}, as well as \cref{it:box-exist} (which determines when the box dimension exists), and \cref{it:rate-form} (describing the covering class of the limit set).
Finally, in \cref{s:examples} we construct examples showing sharpness of certain bounds and completing the proof of \cref{it:fixed}.
We also construct examples of sets of continued fraction expansions with restricted digits (which are invariant for the Gauss map) and prove \cref{it:gauss-invariant}.

\section{An asymptotic formula for covering numbers}\label{s:proof}
In this section, we prove our core result, which is \cref{it:main}.
\subsection{Bounded neighbourhood condition}\label{s:separation}
In the introduction, we assumed that a \cifs{} satisfies the open set condition and the cone condition.
In fact, throughout this section, the only separation assumption we will require is given in \cref{d:bnc} below.

We say that a subset $\mathcal{F}\subset\mathcal{I}^*$ is \emph{mutually incomparable} if $\mtt{i}$ is not a prefix of $\mtt{j}$ for all $\mtt{i},\mtt{j}\in\mathcal{F}$ with $\mtt{i}\neq\mtt{j}$.
\begin{definition}\label{d:bnc}
    We say that the IFS $\{S_i\}_{i\in\mathcal{I}}$ satisfies the \emph{bounded neighbourhood condition} if there exists $M\in\N$ so that for all mutually incomparable $\mathcal{F}\subset\mathcal{I}^*$, for all $x\in X$, and for all $r\in(0,1)$,
    \begin{equation*}
        \#\{\mtt{i}\in\mathcal{F}:\rho(\mtt{i})>r\text{ and }S_{\mtt{i}}(X)\cap B(x,r)\neq\varnothing\}\leq M.
    \end{equation*}
\end{definition}
By a measure argument, it is straightforward to see that the open set condition and the cone condition together imply the bounded neighbourhood condition (see, for instance, \cite[Lemma~2.7]{zbl:0852.28005}).

In some sense, we require the open set condition and the cone condition (via \cite{zbl:0852.28005}) to ensure that
\begin{equation*}
    \dimH\limitset = h\coloneqq\inf\{t \geq 0 : \pressure(t) < 0\}.
\end{equation*}
Under the bounded neighbourhood condition, our proofs hold with the number $h$ in place of $\dimH\limitset$ and do not require any of the results from \cite{zbl:0852.28005}.
We emphasise that the bounded neighbourhood condition is morally very similar to the open set condition and the cone condition; it is likely the case that the results of \cite{zbl:0852.28005} continue to hold only under the assumption of the bounded neighbourhood condition.

\subsection{Regularity of covering numbers}
In this section, we note two standard bounds on covering numbers.

The first bound is an immediate consequence of Ahlfors--David regularity of $\R^d$.
\begin{lemma}\label{l:covering-reg}
    Let $d\in\N$.
    Then there is a constant $A_d\geq 0$ so that for all non-empty bounded sets $E\subset\R^d$ and all $0<\theta\leq 1$,
    \begin{equation*}
        \theta s_E(r^\theta) \leq s_E(r) \leq d - (d-s_E(r^\theta))\theta + \frac{A_d}{\log(1/r)}.
    \end{equation*}
\end{lemma}
\begin{proof}
    On the one hand, $N_r(E)$ increases as $r$ decreases.
    On the other hand, since $\R^d$ is Ahlfors--David $d$-regular, there is a constant $C\geq 1$ (depending on $d$) so that for all $0<r\leq R$, each ball of radius $R$ can be covered by $C(R/r)^{d}$ balls of radius $r$.
    Thus for $0<r<1$, covering balls $B(x,r^\theta)$ by balls of radius $r$,
    \begin{equation*}
        N_{r^\theta}(E) \leq N_{r}(E)\leq C N_{r^{\theta}}(E) \left(\frac{r^\theta}{r}\right)^{d}.
    \end{equation*}
    Writing $A_d = \log(C)$, taking logarithms and rearranging yields the claim.
\end{proof}

The second bound simply uses the definition of the upper box dimension.
\begin{lemma}\label{l:small-r-reg}
    Let $d\in\N$ and $E\subset\R^d$ be non-empty and bounded.
    Then
    \begin{equation*}
        \lim_{r\to 0}\sup_{\theta\in (0,1]}(\theta s_E(r^\theta)-\theta\dimuB E)=0.
    \end{equation*}
    In particular, for all $0\leq h\leq d$,
    \begin{equation*}
        \limsup_{r\to 0}\sup_{\theta\in(0,1]}\left((1-\theta)\cdot h + \theta s_E(r^\theta)\right)=\max\{h,\dimuB E\}.
    \end{equation*}
\end{lemma}
\begin{proof}
    Let $\varepsilon>0$ be arbitrary and let $r_0$ be such that $s_E(r)\leq \dimuB E+\varepsilon$ for all $0<r\leq r_0$.
    Let $M=N_{r_0}(E)$ and let $0<r\leq r_0^{1/\varepsilon}$ be sufficiently small so that
    \begin{equation*}
        \frac{\log M}{\log(1/r)}\leq \varepsilon.
    \end{equation*}
    Then for all $0<\theta\leq 1$, if $\varepsilon\leq \theta\leq 1$, then
    \begin{equation*}
        \theta s_E(r^\theta) \leq \theta (\dimuB E+\varepsilon) \leq \theta \dimuB E+\varepsilon
    \end{equation*}
    and if $0<\theta\leq \varepsilon$, then
    \begin{equation*}
        \theta s_E(r^\theta) = \frac{\log N_{r^\theta}(E)}{\log(1/r)}\leq \frac{\log M}{\log(1/r)}\leq \varepsilon\leq \theta\dimuB E+\varepsilon.
    \end{equation*}
    Thus for all $r$ sufficiently small,
    \begin{equation*}
        \sup_{\theta\in(0,1]} (\theta s_E(r^\theta)-\theta\dimuB E)\leq \varepsilon.
    \end{equation*}

    Conversely, let $r_0>0$ be such that $s_E(r_0)\geq \dimuB E-\varepsilon$ and let $0<r\leq r_0$ be arbitrary.
    Let $\theta_0$ be such that $r^{\theta_0} = r_0$.
    Then
    \begin{equation*}
        \sup_{\theta\in(0,1)} (\theta s_E(r^\theta)-\theta\dimuB E) \geq \theta_0(\dimuB E-\varepsilon) - \theta_0\dimuB E\geq -\varepsilon.
    \end{equation*}
    Since $\varepsilon>0$ was arbitrary, the first claim follows.

    To prove the second claim, by considering $\theta=1$ and $\theta\to 0$ for each $r$, we see that
    \begin{equation*}
        \limsup_{r\to 0}\sup_{\theta\in(0,1]}\left((1-\theta)\cdot h + \theta s_E(r^\theta)\right)\geq \max\{h,\dimuB E\}.
    \end{equation*}
    To obtain the upper bound, let $\varepsilon>0$ be arbitrary.
    Then for $r>0$ sufficiently small, by the first claim,
    \begin{align*}
        \sup_{\theta\in(0,1]}\left((1-\theta)\cdot h + \theta s_E(r^\theta)\right) &\leq \sup_{\theta\in(0,1]}\left((1-\theta)\cdot h + \theta \dimuB E\right) + \varepsilon\\
                             &\leq \max\{h,\dimuB E\}+\varepsilon.
    \end{align*}
    Since $\varepsilon>0$ was arbitrary, the second claim follows.
\end{proof}

\subsection{Properties of the lower box dimension formula}\label{ss:form}
In the statement of \cref{it:main}, an asymptotic formula for $N_r(\limitset)$ is established in terms of the function $\psi(r)$.
In particular, $\psi(r)$ depends implicitly on the precise choice $\fixedpts$.
By taking a higher power, the associated set of fixed points is larger so the associated estimate could in principle be different.
Moreover, one might also consider alternatives to the set of fixed points, such as the \emph{orbit sets}
\begin{equation*}
    \mathcal{O}(x,m)\coloneqq \{S_{\mtt{i}}(x):\mtt{i}\in\mathcal{I}^m\}
\end{equation*}
defined for $x\in X$ and $n\in\N$.
For example, this is the choice made in \cite{zbl:0940.28009}.
Of course, since none of these operations change $N_r(\limitset)$, if \cref{it:main} is to be true, then they also certainly cannot change the asymptotics of the corresponding definition of $\psi(r)$.

In this section, we establish some basic properties of the function $\psi$.
In particular, we will see that there is flexibility in the choice of the set $\fixedpts$, and moreover show that $\psi(r)$ does not depend on the initial level of iteration (up to some error term which is asymptotically $0$).
These results will also be required in our proof of \cref{it:main}.

We now make our definitions precise.
Let $E\subset X$ be an arbitrary bounded non-empty set.
We define
\begin{align*}
    \Psi_E(r,\theta) &= (1-\theta)\dimH\limitset + \frac{\log N_{r^\theta}(E)}{\log(1/r)};\\
    \psi_E(r) &= \liminf_{r\to 0}\max_{\theta\in[0,1]}\Psi_E(r,\theta).
\end{align*}
The maximum in the definition of $\psi_E(r)$ is attained because the covering number $N_r$ is defined in terms of open balls, so the map $\theta\mapsto N_{r^\theta}(E)$ is upper semi-continuous.
\begin{definition}
    Let $\{S_i\}_{i\in\mathcal{I}}$ be a \cifs{} on a compact set $X$.
    We say a set $E \subset \bigcup_{i \in \mathcal{I}} S_i(X)$ is a \emph{discrete approximation} of $\{S_i\}_{i\in\mathcal{I}}$ if there is a number $k\in\N$ so that $1\leq \# (E\cap S_i(x))\leq k$ for all $i\in\mathcal{I}$.
\end{definition}
For instance, the set of fixed points $\fixedpts$ and the orbit sets $\mathcal{O}(x,1)$ for $x\in X$ are discrete approximations.
\begin{lemma}\label{l:fixed-choice}
    Let $\{S_i\}_{i\in\mathcal{I}}$ be a \cifs{} on a compact set $X$.
    Suppose $E_1$ and $E_2$ are discrete approximations of $\{S_i\}_{i\in\mathcal{I}}$.
    Then $N_r(E_1) \approx N_r(E_2)$, uniformly for $r \in (0,1)$.
    In particular,
    \begin{equation*}
        \lim_{r\to 0}\bigl(\psi_{E_1}(r)-\psi_{E_2}(r)\bigr)=0.
    \end{equation*}
\end{lemma}
\begin{proof}
    This holds by the same argument as the proof of \cite[Proposition~2.9]{zbl:0940.28009}, which depends on the bounded neighbourhood condition.
\end{proof}
Next, we consider higher iterates of the \cifs.
With \cref{l:fixed-choice} in mind, for convenience we may restrict our attention to the orbit sets $\mathcal{O}(x_0,m)$, where $m\in\N$, and $x_0\in X$ is chosen so that $S_i(x_0)=x_0$ for some $i\in\mathcal{I}$.
Write $F_m\coloneqq \mathcal{O}(x_0,m)$.
With $x_0$ chosen in this way,
\begin{equation}\label{e:F-invariance}
    F_m=\bigcup_{\mtt{i}\in\mathcal{I}^{m-1}}S_{\mtt{i}}(F_1) \subset \limitset.
\end{equation}
For each $m\in\N$, as shorthand, we also write $\Psi_m\coloneqq \Psi_{F_m}$, $\psi_m\coloneqq \psi_{F_m}$, and for $0<r<1$
\begin{equation*}
    s_m(r)\coloneqq \frac{\log N_{r}(\fixedpts_m)}{\log(1/r)}.
\end{equation*}
Next, for each $\mtt{i}\in\mathcal{I}^*$ with $\rho(\mtt{i})>r$, let $\theta_{\mtt{i}}(r)\in(0,1]$ be given by
\begin{equation*}
    \theta_{\mtt{i}}(r) \coloneqq 1-\frac{\log \rho(\mtt{i})}{\log r}.
\end{equation*}
Equivalently, $\theta_{\mtt{i}}(r)$ is chosen so that
\begin{equation}\label{e:lambda-choice}
    N_{r\cdot \rho(\mtt{i})^{-1}}(\fixedpts_1) =  \rho(\mtt{i})^{\dimH\limitset}\cdot\left(\frac{1}{r}\right)^{\Psi_1\left(r,\theta_{\mtt{i}}(r)\right)}.
\end{equation}
For $m\in\N$ and $0<r<1$, we also write
\begin{equation*}
    E_m(r)\coloneqq \bigcup_{\substack{\mtt{i}\in\mathcal{I}^m\\\rho(\mtt{i})\leq r}}S_{\mtt{i}}(\limitset).
\end{equation*}
Since $E_m(r)$ is contained in the $(r\cdot\diam X)$-neighbourhood of $\fixedpts$,
\begin{equation}\label{e:C-choice}
    N_r(E_m(r))\lesssim N_r(\fixedpts).
\end{equation}
We now establish invariance under higher iterates.
\begin{lemma}\label{l:higher}
    Let $\{S_i\}_{i\in\mathcal{I}}$ be a \cifs.
    Then for each $m\in\N$,
    \begin{equation*}
        \lim_{r\to 0}\bigl(\psi_1(r)-\psi_m(r)\bigr)=0.
    \end{equation*}
\end{lemma}
\begin{proof}
    It is clear that $\psi_1\leq \psi_2\leq\cdots\leq \psi_m$ for $m\in\N$.
    Thus it suffices to prove that for all $\varepsilon>0$ and all $m\in\N$ sufficiently large,
    \begin{equation}\label{e:limsup-version}
        \limsup_{r\to 0}\bigl(\psi_m(r)-\psi_1(r)\bigr)\leq 2\varepsilon.
    \end{equation}

    Fix $\varepsilon>0$.
    By the definition of the pressure, let $m\in\N$ be sufficiently large so that
    \begin{equation}\label{e:pressure-sum}
        \sum_{\mtt{i}\in\mathcal{I}^{m-1}} \rho(\mtt{i})^{\dimH\limitset+\varepsilon}<\infty.
    \end{equation}
    We first show that
    \begin{equation}\label{e:sn-bound}
        \limsup_{r\to 0}\bigl(s_m(r)-\psi_1(r)\bigr)\leq \varepsilon.
    \end{equation}
    Let $\varepsilon>0$ be arbitrary.
    For each $r\in(0,1)$, recalling \cref{e:F-invariance} and \cref{e:C-choice},
    \begin{align*}
        N_r(F_n) &= N_r\left(\bigcup_{\mtt{i}\in\mathcal{I}^{m-1}} S_{\mtt{i}}(F_1)\right)\\
                 &\lesssim N_r(F_1) + \sum_{\substack{\mtt{i}\in\mathcal{I}^{m-1}\\\rho(\mtt{i})>r}}N_{r\cdot \rho(\mtt{i})^{-1}}(\fixedpts_1)\\
                 &= N_r(F_1) + \sum_{\substack{\mtt{i}\in\mathcal{I}^{m-1}\\\rho(\mtt{i})>r}} \rho(\mtt{i})^{\dimH\limitset}\cdot\left(\frac{1}{r}\right)^{\Psi_1\left(r,\theta_{\mtt{i}}(r)\right)}\\
                 &\leq N_r(F_1) + \left(\frac{1}{r}\right)^{\psi_1(r)+\varepsilon}\sum_{\mtt{i}\in\mathcal{I}^{m-1}} \rho(\mtt{i})^{\dimH\limitset+\varepsilon}\\
                 &\lesssim \left(\frac{1}{r}\right)^{\psi_1(r)+\varepsilon}.
    \end{align*}
    In the last line, we used \cref{e:pressure-sum} and the observation that $s_1(r)=\Psi_1(r,1)\leq \psi_1(r)$.
    The above calculation is equivalent to the fact that there is a constant $C>0$ so that
    \begin{equation*}
        s_m(r) \leq \psi_1(r)+\frac{C}{\log(1/r)}+\varepsilon.
    \end{equation*}
    Thus \cref{e:sn-bound} follows.

    We now establish \cref{e:limsup-version}.
    Let $r\in(0,1)$ be small, let $\theta$ be chosen so that $\psi_m(r)=\Psi_m(r,\theta)$, and let $\kappa$ be chosen so that $\psi_1(r)=\Psi_1(r,\kappa)$.
    By the definition of $\kappa$ and using \cref{e:sn-bound}, for all $r$ sufficiently small,
    \begin{align*}
        \psi_m(r) &= (1-\theta) \dimH\limitset + \theta s_m(r^\theta)\\
                  &\leq (1-\theta)\dimH\limitset + \theta \bigl((1-\kappa)\dimH\limitset+\kappa s_1(r^{\theta\kappa})+2\varepsilon\bigr)\\
                  &= \Psi_1(r,\theta\kappa) + 2\theta\varepsilon\\
                  &\leq \psi_1(r) +2\varepsilon.
    \end{align*}
    Thus \cref{e:limsup-version} follows, and therefore the desired result holds.
\end{proof}

\subsection{Proof of the asymptotic formula}
In this section, we establish our main asymptotic formula for $N_r(\limitset)$, as stated in \cref{it:main}.

First, for each $m\in\N$, define
\begin{equation*}
    \tau_m(r)
    \coloneqq\sum_{\substack{\mtt{i}\in(\mathcal{I}^m)^*\\\rho(\mtt{i})>r}}N_{r\cdot \rho(\mtt{i})^{-1}}(\fixedpts_m)
    =\sum_{\substack{\mtt{i}\in(\mathcal{I}^m)^*\\\rho(\mtt{i})>r}} \rho(\mtt{i})^{\dimH\limitset}\left(\frac{1}{r}\right)^{\Psi_m\left(r,\theta_{\mtt{i}}(r)\right)}.
\end{equation*}
We first reduce the question of bounding the covering numbers $N_r(\limitset)$ to the question of bounding the symbolic counts $\tau_m(r)$.
\begin{lemma}\label{l:tau-bound}
    Let $\{S_i\}_{i\in\mathcal{I}}$ be a \cifs{} with attractor $\limitset$ and fixed points $\fixedpts$.
    Then for all $\varepsilon>0$ and $m\in\N$ sufficiently large, there exists $C\geq 1$ so that for all $0<r<1$,
    \begin{equation*}
        C^{-1} r^{\varepsilon}\tau_m(r)\leq N_r(\limitset)\leq C r^{-\varepsilon}\tau_m(r).
    \end{equation*}
\end{lemma}
\begin{proof}
    The result will follow from the following key estimate: there exists a constant $C_0\geq 1$ so that for all $m\in\N$ and $0<r<1$,
    \begin{equation}\label{e:symbolic-transfer}
        C_0^{-1}\cdot N_r(\limitset) \leq N_r(\fixedpts_m) + \sum_{\substack{\mtt{i}\in\mathcal{I}^m\\\rho(\mtt{i})>r}}N_{r\cdot \rho(\mtt{i}^{-1})}(\limitset)\leq C_0\cdot N_r(\limitset)
    \end{equation}
    Critically, the constant $C_0$ is independent of $m$ and $r$.
    Suppose this estimate holds.
    Recalling that $\xi=\sup_{i\in\mathcal{I}}\rho(i)\in(0,1)$, take $k\in\N$ be minimal so that $\xi^{km}\leq r$.
    Then applying \cref{e:symbolic-transfer} $k$ times gives,
    \begin{equation*}
        N_r(\limitset) \leq \sum_{\substack{\mtt{i}\in(\mathcal{I}^m)^*\\\rho(\mtt{i})>r\text{ and }|\mtt{i}| < k}}C_0^{|\mtt{i}|+1}\cdot N_{r\cdot \rho(\mtt{i})^{-1}}(\fixedpts_m)\leq C_0^k \tau_m(r).
    \end{equation*}
    The lower bound $C_0^{-k}\tau_m(r)\leq N_r(\limitset)$ also holds by the same argument.
    Then, given $\varepsilon>0$, let $m$ be sufficiently large so that $\frac{\log C_0}{m\log(1/\xi)}\leq\varepsilon$ and, since $\xi^{(k-1)m}>r$,
    \begin{equation*}
        C_0^{-1} r^\varepsilon\tau_m(r)\leq C_0^{-k}\tau_m(r)\leq N_r(\limitset)\leq C_0^k\tau_m(r)\leq  C_0 r^{-\varepsilon}\tau_m(r)
    \end{equation*}
    which is the desired result.

    It remains to establish \cref{e:symbolic-transfer}.
    By the invariance property of $\limitset$,
    \begin{equation*}
        \limitset = E_m(r)\cup \bigcup_{\substack{i\in\mathcal{I}^m\\\rho(\mtt{i})>r}}S_{\mtt{i}}(\limitset).
    \end{equation*}
    Thus by \cref{e:C-choice}, since $N_r(S_{\mtt{i}}(\limitset))\approx N_{r\cdot \rho(\mtt{i})^{-1}}(\limitset)$,
    \begin{equation*}
        N_r(\limitset) \lesssim N_r(\fixedpts_m) + \sum_{\substack{i\in\mathcal{I}^m\\\rho(\mtt{i}) > r}} N_{r\cdot \rho(\mtt{i})^{-1}}(\limitset).
    \end{equation*}

    We now obtain the other bound.
    First, by the bounded neighbourhood condition,
    \begin{equation*}
        \sum_{\substack{\mtt{i}\in\mathcal{I}^m\\ \rho(\mtt{i})>r}} N_{r}(S_{\mtt{i}}(\limitset))\leq M\cdot N_r\left(\bigcup_{\substack{\mtt{i}\in\mathcal{I}^m\\ \rho(\mtt{i})>r}}S_{\mtt{i}}(\limitset)\right)\leq M\cdot N_r(\limitset).
    \end{equation*}
    Therefore since $N_r(\fixedpts_m)\leq N_r(\limitset)$,
    \begin{equation*}
        N_r(\fixedpts_m) + \sum_{\substack{\mtt{i}\in\mathcal{I}^m\\\rho(\mtt{i})>r}}N_{r\cdot \rho(\mtt{i})^{-1}}(\limitset)\lesssim (M+1)\cdot N_r(\limitset).
    \end{equation*}
    Thus \cref{e:symbolic-transfer} follows.
\end{proof}
Next, we require the standard observation for finite iterated function systems that the Hausdorff dimension is realised uniformly over all scales.
\begin{lemma}\label{l:saturation}
    Let $\{S_i\}_{i\in\mathcal{I}}$ be a \cifs{} with limit set $\limitset$, and suppose $\mathcal{I}$ is a finite index set.
    Let $\rho_{\min}=\min\{\rho(i):i\in\mathcal{I}\}>0$.
    Then for all $s<\dimH\limitset$, with
    \begin{equation*}
        \mathcal{I}(r)\coloneqq \{\mtt{i}\in\mathcal{I}^*: r<\rho(\mtt{i})\leq r K \rho_{\min}^{-1}\},
    \end{equation*}
    we have $\#\mathcal{I}(r)\gtrsim r^{-s}$.
\end{lemma}
\begin{proof}
    Since $\mathcal{I}$ is finite, $\dimB\limitset = \dimH\limitset$.
    Moreover, by \cref{e:rho-approx-mul}, every infinite word $\gamma\in\mathcal{I}^{\N}$ has at least one prefix in $\mathcal{I}(r)$, so $\{S_{\mtt{i}}(X):\mtt{i}\in\mathcal{I}(r)\}$ is a cover for $\limitset$.

    Now fix $s<\dimB\limitset$ and let $r>0$ be small.
    Get a pairwise-disjoint family of balls $\{B(x_i,r)\}_{i=1}^N$ with $N\gtrsim r^{-s}$ and each $x_i\in X$.
    Thus by the bounded neighbourhood condition, $\#\mathcal{I}(r)\geq M^{-1} N\gtrsim r^{-s}$, as claimed.
\end{proof}
We also need a simple continuity property of the function $\Psi(r,\theta)$ in $\theta$.
\begin{lemma}\label{l:psi-cont}
    There is a constant $A_d\geq 0$ so that for all $m\in\N$, $\theta_1,\theta_2\in[0,1]$ and $0<r<1$,
    \begin{equation*}
        \left\lvert\Psi_m(r,\theta_1)-\Psi_m(r,\theta_2)\right\rvert \leq 2d |\theta_1-\theta_2|+\frac{A_d}{\log(1/r)}.
    \end{equation*}
\end{lemma}
\begin{proof}
    Let $m\in\N$ and $0<r<1$.
    Without loss of generality, we may fix $0\leq \theta_1\leq\theta_2\leq 1$.
    Then by \cref{l:covering-reg},
    \begin{align*}
        \left\lvert\Psi_m(r,\theta_1)-\Psi_m(r,\theta_2)\right\rvert &\leq (\theta_2-\theta_1)\dimH\limitset + \left\lvert \frac{\log N_{r^{\theta_2}}(F_m)}{\log(1/r)}-\frac{\log N_{r^{\theta_1}}(F_m)}{\log(1/r)}\right\rvert\\
                                                                     &\leq (\theta_2-\theta_1) d + \frac{A_d}{\log(1/r)} + (\theta_2-\theta_1) d
    \end{align*}
    as claimed.
\end{proof}
We now have all of the tools required to prove our main formula, which we restate below for the convenience of the reader.
\begin{restatement}{it:main}
    Let $\limitset$ be the limit set of a \cifs{} on $\R^d$ with fixed points $F$ and associated function $\psi$.
    Then
    \begin{equation*}
        \lim_{r\to 0}\left(\frac{\log N_r(\limitset)}{\log(1/r)}-\psi(r)\right)=0.
    \end{equation*}
    In particular,
    \begin{equation*}
        \dimlB\limitset = \liminf_{r\to 0}\psi(r).
    \end{equation*}
\end{restatement}
\begin{proof}
    By \cref{l:fixed-choice,l:higher}, it suffices to show that for all $\varepsilon>0$, there exists $m\in\N$ so that for all $r>0$ sufficiently small,
    \begin{equation*}
        -6\varepsilon\leq \frac{\log N_r(\limitset)}{\log(1/r)}-\psi_m(r)\leq 2\varepsilon.
    \end{equation*}

    We begin with the upper bound.
    First, let $m\in\N$ be sufficiently large so that
    \begin{equation}\label{e:above-P}
        \sum_{\mtt{i}\in\mathcal{I}^m}\rho(\mtt{i})^{\dimH\limitset+\varepsilon}\eqqcolon \vartheta<1,
    \end{equation}
    and moreover for all $r>0$ sufficiently small
    \begin{equation*}
        \frac{\log N_r(\limitset)}{\log(1/r)}\leq \frac{\log \tau_m(r)}{\log(1/r)}+\varepsilon.
    \end{equation*}
    The second choice is possible by \cref{l:tau-bound}.
    The choice \cref{e:above-P} implies by sub-multiplicativity of $\rho$ that
    \begin{align*}
        \sum_{\mtt{i}\in(\mathcal{I}^m)^*} \rho(\mtt{i})^{\dimH\limitset+\varepsilon} &= \sum_{k=0}^\infty \sum_{\mtt{i}_1\in\mathcal{I}^m}\cdots \sum_{\mtt{i}_k\in\mathcal{I}^m}\rho(\mtt{i}_1\cdots\mtt{i}_k)^{\dimH\limitset+\varepsilon}\\
                                                                                      &\leq \sum_{k=0}^\infty \left(\sum_{\mtt{i}\in\mathcal{I}^m} \rho(\mtt{i})^{\dimH\limitset+\varepsilon}\right)^k\\
                                                                                      & =\sum_{k=0}^\infty\vartheta^k<\infty.
    \end{align*}
    Therefore
    \begin{align*}
        \tau_m(r)&=\sum_{\substack{\mtt{i}\in(\mathcal{I}^m)^*\\ \rho(\mtt{i})>r}} \rho(\mtt{i})^{\dimH\limitset}\cdot\left(\frac{1}{r}\right)^{\Psi_m(r,\theta_{\mtt{i}}(r))}\\
                 &\leq \left(\frac{1}{r}\right)^{\psi_m(r)+\varepsilon}\sum_{\substack{\mtt{i}\in(\mathcal{I}^m)^*\\ \rho(\mtt{i})>r}} \rho(\mtt{i})^{\dimH\limitset+\varepsilon}\\
                 &\lesssim_\varepsilon \left(\frac{1}{r}\right)^{\psi_m(r)+\varepsilon}.
    \end{align*}

    We now establish the lower bound.
    Heuristically, the upper bound proved above is sharp if the sum of the $\rho(\mtt{i})^{\dimH\limitset}$ is realised uniformly over all scales simultaneously so that the error which results from bounding $\Psi_m(r,\cdot)$ by the supremum $\psi_m(r)$ is small.
    This follows by approximation via finite subsystems as a direct consequence of \cref{l:saturation}.

    Let $m\in\N$ be sufficiently large so that there is a finite subset $\mathcal{F}\subset\mathcal{I}^m$ such that the \cifs{} $\{S_{\mtt{i}}\}_{\mtt{i}\in\mathcal{F}}$ has limit set $\limitset_\varepsilon$ with
    \begin{equation*}
        s\coloneqq \dimH\limitset_\varepsilon>\dimH\limitset-\varepsilon,
    \end{equation*}
    and moreover for all $r>0$ sufficiently small
    \begin{equation*}
        \frac{\log N_r(\limitset)}{\log(1/r)}\geq \frac{\log \tau_m(r)}{\log(1/r)}-\varepsilon.
    \end{equation*}
    That the first choice is possible follows from \cite[Theorem~3.15]{zbl:0852.28005}, and again the second choice is possible by \cref{l:tau-bound}.
    Next, fix the constant $A_d$ from \cref{l:psi-cont}, let $r$ be sufficiently small so that $A_d/\log(1/r)\leq\varepsilon$, and fix a partition $0=\kappa_0<\kappa_1<\cdots<\kappa_\ell = 1$ such that $\kappa_i-\kappa_{i-1}\leq (2d)^{-1}\varepsilon$ for all $i=1,\ldots,\ell$.
    Now for each $i=1,\ldots,\ell$ and $r\in(0,1)$, set
    \begin{align*}
        \mathcal{F}^*_i(r) &= \bigl\{\mtt{i}\in\mathcal{F}^*: \rho(\mtt{i})>r;\, \kappa_{i-1} \leq 1 - \theta_{\mtt{i}}(r) < \kappa_i\bigr\}\\
                           &= \bigl\{\mtt{i}\in\mathcal{F}^*: r^{\kappa_{i}} < \rho(\mtt{i})\leq r^{\kappa_{i-1}}\bigr\}.
    \end{align*}
    Note that if $\mtt{i}\in\mathcal{F}^*_i(r)$, by \cref{l:psi-cont} and the choice of the $\kappa_i$, for all $r$ sufficiently small,
    \begin{equation}\label{e:kappa-approx}
        \Psi_m(r,\theta_{\mtt{i}}(r)) \geq \Psi_m(r,1-\kappa_i)-2\varepsilon.
    \end{equation}
    Now, let $r$ moreover be sufficiently small so that for all $i=1,\ldots,\ell$,
    \begin{equation*}
        K\rho_{\mathcal{F}}^{-1} r^{\kappa_i}\leq r^{\kappa_{i-1}}\quad\text{where}\quad \rho_{\mathcal{F}}\coloneqq \min\{\rho(\mtt{i}):\mtt{i}\in\mathcal{F}\}.
    \end{equation*}
    Thus for such $r$ and all $i=1,\ldots,\ell$, by \cref{l:saturation},
    \begin{align*}
        \#\mathcal{F}^*_i(r) &\geq \#\{\mtt{i}\in\mathcal{F}^*:r^{\kappa_i} < \rho(\mtt{i}) \leq K \rho_{\mathcal{F}}^{-1} r^{\kappa_i}\}\gtrsim_\varepsilon r^{-\kappa_i (\dimH\limitset-\varepsilon)},
    \end{align*}
    hence
    \begin{equation*}
        \sum_{\mtt{i}\in \mathcal{F}^*_i(r)} \rho(\mtt{i})^{\dimH\limitset-\varepsilon}\geq \#\mathcal{F}^*_i(r) r^{\kappa_i(\dimH\limitset-\varepsilon)}\gtrsim_\varepsilon 1.
    \end{equation*}
    Thus for all $r$ sufficiently small,
    \begin{align*}
        \tau_m(r)&=\sum_{\substack{\mtt{i}\in(\mathcal{I}^m)^*\\ \rho(\mtt{i})>r}} \rho(\mtt{i})^{\dimH\limitset}\cdot\left(\frac{1}{r}\right)^{\Psi_m(r,\theta_{\mtt{i}}(r))}\\
                 &\geq \sum_{i=1}^\ell \left(\frac{1}{r}\right)^{\Psi_m(r,1-\kappa_i)-3\varepsilon} \sum_{\mtt{i}\in \mathcal{F}^*_i(r)} \rho(\mtt{i})^{\dimH\limitset-\varepsilon}\\
                 &\gtrsim_{\varepsilon} \max_{i=1,\ldots,\ell}\left(\frac{1}{r}\right)^{\Psi_m(r,1-\kappa_i)-3\varepsilon}\\
                 &\geq \left(\frac{1}{r}\right)^{\psi_m(r)-5\varepsilon}.
    \end{align*}
    In the last line, we again applied \cref{l:psi-cont} using the observation that $\{\kappa_1,\ldots,\kappa_l\}$ is $(2d)^{-1}\varepsilon$-dense in $[0,1]$.
    Since $\varepsilon>0$ was arbitrary, the desired result follows.
\end{proof}

In \cite[Theorem~3.1]{zbl:0852.28005} it was shown that the packing dimension, which is the same as the modified upper box dimension, always coincides with upper box dimension for the attractor of a \cifs.
Of course, the analogous result holds for modified lower box dimension, which is defined by
\begin{equation*}
    \dimlMB K = \inf\Big\{\sup_i \dimlB K_i : K \subset \bigcup_{i=1}^{\infty} K_i \Big\}.
\end{equation*}
This is a standard consequence of the Baire category theorem (see, e.g., \cite[Proposition~2.8]{zbl:1285.28011}).
\begin{proposition}
    Let $\limitset$ be the limit set of a \cifs.
    Then $\dimlMB \limitset = \dimlB \limitset$.
\end{proposition}

\begin{remark}
    Given a countably infinite CIFS, it may be of interest to consider the infimum of dimensions of limit sets of cofinite subsystems.
    In our case, for instance, the lower box dimension converges to the same formula as in \cref{it:main}, except with the \emph{finiteness parameter} $\inf\{t \geq 0 : \pressure(t) < \infty\}$ in place of the Hausdorff dimension in the definition of $\Psi$.
    Similar statements can also be made about the asymptotics of the covering numbers of the cofinite subsystems.
\end{remark}

\section{Consequences of the asymptotic formula}\label{s:consequences}
In this section, we obtain consequences of the asymptotic formula stated in \cref{it:main}.
\subsection{Classifying existence of the box dimension}
Using our formula for the lower box dimension stated in \cref{it:main}, we obtain bounds on the lower box dimension in terms of $\dimH\limitset$, $\dimlB\fixedpts$, $\dimuB\fixedpts$, and the ambient dimension $d$, without any other information concerning the set $\fixedpts$.
Recalling the general bounds from \cref{e:trivialbounds}, this implies the first half of \cref{it:fixed}.
\begin{corollary}\label{c:fixed-bound}
    Let $\limitset$ be the limit set of a \cifs{} on $\R^d$.
    Then if $\dimuB\fixedpts>\dimH\limitset$,
    \begin{equation*}
        \dimlB \limitset\leq \dimH\limitset + \frac{(\dimuB\fixedpts - \dimH\limitset) (d - \dimH\limitset) \dimlB\fixedpts}{d\dimuB\fixedpts-\dimH\limitset \cdot\dimlB\fixedpts}.
    \end{equation*}
\end{corollary}
\begin{proof}
    For notational simplicity, write $s=\dimlB\fixedpts$, $t=\dimuB\fixedpts$, and $h=\dimH\limitset$.
    Note that $h<t\leq d$, so we may set
    \begin{equation*}
        \theta_d \coloneqq \frac{(d-h) t}{d\cdot t-h\cdot s}\leq 1.
    \end{equation*}
    Let $(r_n)_{n=1}^\infty$ be a sequence converging to zero such that
    \begin{equation*}
        \dimlB\fixedpts = \liminf_{r\to 0}s_F(r)=\lim_{n\to\infty}s_F(r_n^{\theta_d}).
    \end{equation*}

    Let $\varepsilon>0$ and let $n$ be sufficiently large so that $s_F(r_n^{\theta_d})\leq s +\varepsilon$ and $\theta s_F(r_n^\theta) \leq \theta\dimuB\fixedpts+\varepsilon$ for all $0<\theta\leq 1$ by \cref{l:small-r-reg}, and $A_d \leq \varepsilon \log(1/r_n^\theta)$ where $A_d$ is the constant from \cref{l:covering-reg}.
    By \cref{it:main} and the definition of $\theta_d$, it suffices to show that for all $0<\theta\leq 1$,
    \begin{equation*}
        (1-\theta)\cdot h + \theta s_F(r_n^\theta) \leq d - (d-s)\theta_d+ 2\varepsilon.
    \end{equation*}
    We consider three cases depending on the value of $\theta$.
    \begin{enumerate}[nl]
        \item $\theta_d\leq \theta\leq 1$.
            Then by \cref{l:covering-reg},
            \begin{align*}
                s_F(r_n^\theta) &\leq d - (d - s_F(r_n^{\theta_d}))\frac{\theta_d}{\theta} + \varepsilon \\
                              &\leq d - (d-s)\theta_d + 2\varepsilon.
            \end{align*}
            Thus since $t>h$,
            \begin{align*}
                (1-\theta) h + \theta s_F(r_n^\theta) &\leq \max\{h, d - (d-s)\theta_d\}+2\varepsilon\\
                                                    &=d - (d-s)\theta_d+2\varepsilon.
            \end{align*}
        \item $\theta_d\cdot s/t < \theta < \theta_d$.
            Then by \cref{l:covering-reg},
            \begin{align*}
                (1-\theta) h + \theta s_F(r_n^\theta) &\leq (1-\theta) h +\theta_d s_F(r_n^{\theta_d})\\
                                                    &\leq (1-\theta) h + \theta_d s + \theta_d\varepsilon\\
                                                    &\leq d - (d-s)\theta_d +\varepsilon.
            \end{align*}
            Here, the final inequality is equivalent to the lower bound on $\theta$.
        \item $0<\theta\leq \theta_d \cdot s/t$.
            Then since $t>h$,
            \begin{align*}
                (1-\theta) h + \theta s_F(r_n^\theta) &\leq h +\theta(t-h)+\varepsilon\\
                                                    &\leq h + \theta_d \left(s - \frac{hs}{t}\right)+\varepsilon\\
                                                    & = d-(d-s)\theta_d +\varepsilon.
            \end{align*}
    \end{enumerate}
    This covers all the possible values of $\theta$, as required.
\end{proof}
From this bound, it is straightforward to deduce our main classification result on the existence of the box dimension.
\begin{restatement}{it:box-exist}
    Let $\limitset$ be the limit set of a \cifs.
    Then $\dimlB\limitset=\dimuB\limitset$ if and only if
    \begin{equation*}
        \dimuB\fixedpts\leq \max \{\dimH \limitset, \dimlB \fixedpts\}.
    \end{equation*}
\end{restatement}
\begin{proof}
    First, suppose
    \begin{equation*}
        \dimuB\fixedpts\leq \max \{\dimH \limitset, \dimlB \fixedpts\}.
    \end{equation*}
    Then
    \begin{align*}
        \max\{\dimH\limitset,\dimlB\fixedpts\}&\leq\dimlB\limitset\\
                                              &\leq \dimuB\limitset\\
                                              &=\max\{\dimH\limitset,\dimuB\fixedpts\}\\
                                              &\leq\max\{\dimH\limitset,\dimlB\fixedpts\}
    \end{align*}
    so in fact equality holds, as claimed.

    Conversely, suppose
    \begin{equation*}
        \dimuB\fixedpts>\max \{\dimH \limitset, \dimlB \fixedpts\}.
    \end{equation*}
    Since $\dimuB\limitset=\max\{\dimH\limitset,\dimuB\fixedpts\}$, this implies that
    \begin{equation*}
        \dimuB\limitset=\dimuB\fixedpts\qquad\text{and}\qquad 1-\frac{\dimlB\fixedpts}{\dimuB\fixedpts} > 0.
    \end{equation*}
    Thus by \cref{c:fixed-bound} (or, more precisely, the limiting bound as explained in \cref{r:inf-bound}),
    \begin{align*}
        \dimlB\limitset &\leq \dimlB\fixedpts+\left(1-\frac{\dimlB\fixedpts}{\dimuB\fixedpts}\right)\dimH\limitset\\
                        &<\dimlB\fixedpts + \left(1-\frac{\dimlB\fixedpts}{\dimuB\fixedpts}\right)\dimuB\fixedpts\\
                        &= \dimuB\fixedpts\\
                        &= \dimuB\limitset
    \end{align*}
    as claimed.
\end{proof}

\subsection{Some preliminaries on the covering class}\label{ss:cov-class}
In this section, we provide an introduction to the covering class and in particular prove \cref{p:covering-class}.
We recall the various definitions from \cref{ss:rate-intro}.
In the following lemma we note an equivalent integrated version of the definition of $\mathcal{G}(\lambda,\alpha)$.
This is a consequence of the mean value theorem for one-sided derivatives of continuous functions; for a proof, see for instance \cite[Lemma~3.2]{zbl:1509.28005}.
\begin{lemma}\label{l:scale-bound}
    Let $0\leq\lambda\leq\alpha\leq d$ and let $g \colon \R \to\R$.
    Then $g\in\mathcal{G}(\lambda,\alpha)$ if and only if for all $x_0\in\R$ and $x>0$,
    \begin{equation*}
        \lambda-(\lambda-g(x_0))\exp(-x)\leq g(x_0+x)\leq \alpha-(\alpha-g(x_0))\exp(-x).
    \end{equation*}
\end{lemma}

First, we show that the class $\mathcal{G}(\lambda,\alpha)$ is closed under infima and suprema.
\begin{proposition}\label{p:relative-compact}
    Let $0\leq\lambda\leq\alpha$.
    Then every sequence $(g_n)_{n=1}^\infty\subset\mathcal{G}(\lambda,\alpha)$ has a subsequence which converges uniformly on compact sets to a function $g\in\mathcal{G}(\lambda,\alpha)$.

    Moreover, if $g$ is the pointwise infimum or supremum of a family of functions $g_j\in\mathcal{G}(\lambda,\alpha)$, then $g\in\mathcal{G}(\lambda,\alpha)$.
\end{proposition}
\begin{proof}
    Firstly, the family $\mathcal{G}(\lambda,\alpha)$ is uniformly bounded and uniformly equicontinuous since it is a subset of the set of Lipschitz functions with constant $\alpha-\lambda$ taking values in the interval $[\lambda,\alpha]$.
    Thus by the Arzelà--Ascoli theorem, $g_n$ has a subsequence which converges to a function $g$ uniformly on every compact subset of $\R$.

    Next, we show that $g\in\mathcal{G}(\lambda,\alpha)$.
    It is clear that $g$ takes values in $[\lambda,\alpha]$.
    Let $x_0\in\R$, $x>0$ and $\varepsilon>0$ be arbitrary.
    Then for an infinite sequence of $n$ we have $\sup_{y \in [x_0 - x, x_0 + x]} |g_n(y) - g(y)| \leq \varepsilon$, so by \cref{l:scale-bound} applied to the function $g_n$,
    \begin{align*}
        g(x_0 + x) &\leq g_n(x_0+x)+\varepsilon\\
                   &\leq \alpha - (\alpha - g_n(x_0))\exp(-x) + \varepsilon\\
                   &\leq \alpha - (\alpha - g(x_0))\exp(-x) + \varepsilon(1+\exp(-x)).
    \end{align*}
    Since $\varepsilon>0$ was arbitrary,
    \begin{equation*}
        g(x_0+x)\leq \alpha - (\alpha - g(x_0))\exp(-x).
    \end{equation*}
    The lower bound with $\lambda$ in place of $\alpha$ is identical.
    Thus by \cref{l:scale-bound}, $g\in\mathcal{G}(\lambda,\alpha)$.

    Now let $\{g_j:j\in J\}\subset\mathcal{G}(\lambda,\alpha)$ be an arbitrary family of functions with supremum $g$.
    Since $\mathcal{G}(\lambda,\alpha)$ is uniformly bounded and uniformly equicontinuous, there exists a sequence $J_1\subset J_2\subset\cdots\subset J$ of finite subsets such that
    \begin{equation*}
        g = \lim_{n\to\infty} g_n\quad\text{where}\quad g_n = \max\{g_j:j\in J_n\}.
    \end{equation*}
    Since the sequence $g_n$ is monotonic, it suffices to verify that $g_n\in\mathcal{G}(\lambda,\alpha)$ for each $n\in\N$.
    To check this, it suffices to check that if $f_1,f_2\in\mathcal{G}(\lambda,\alpha)$, then $f\coloneqq \max\{f_1,f_2\}\in\mathcal{G}(\lambda,\alpha)$.
    Let $x\in\R$ be arbitrary.
    Since $f_1$ and $f_2$ are continuous, if $f_1(x)<f_2(x)$, then $\diniu+f(x)=\diniu+f_2(x)$ and $f(x)=f_2(x)$.
    The analogous statement holds if $f_2(x)<f_1(x)$.
    Otherwise if $f_1(x)=f_2(x)$, then $\diniu+ f(x)=\max\{\diniu+ f_1(x),\diniu+f_2(x)\}$.
    In either case it is clear that $\diniu+ f(x)\in[\lambda-f(x),\alpha-f(x)]$, as required.

    The case for the infimum of a family of functions is identical.
\end{proof}
Next, we show that establishing an approximate form of the inequalities in \cref{l:scale-bound} suffices to show asymptotic equivalence to a function in $\mathcal{G}(\lambda,\alpha)$.
\begin{lemma}\label{l:approx-bound}
    Let $e\colon\R\to\R$ be any function with $\lim_{x\to\infty}e(x)=0$ and let $z\in\R$.
    Suppose $0\leq\lambda\leq\alpha$ and let $f\colon\R\to\R$ be any function such that for all $x_0\geq z$ and $x\geq 0$,
    \begin{align*}
        \lambda-(\lambda - f(x_0) + e(x_0))\exp(-x) &\leq f(x_0+x)\\
                                                     &\leq \alpha - (\alpha- f(x_0))\exp(-x) + e(x_0+x).
    \end{align*}
    Then there exists $g\in\mathcal{G}(\lambda,\alpha)$ such that $f\eqinf g$.
\end{lemma}
\begin{proof}
    For $y\geq z$, observe that
    \begin{equation*}
        \lambda(1-\exp(z-y)) + (f(z)-e(z))\exp(z-y)\leq f(y)
    \end{equation*}
    and similarly
    \begin{equation*}
        f(y)\leq \alpha(1-\exp(z-y))+f(z)\exp(z-y)+e(y).
    \end{equation*}
    In particular,
    \begin{equation}\label{e:crude-bound}
        \lambda\leq \liminf_{y\to\infty}f(y)\leq \limsup_{y\to\infty}f(y)\leq \alpha.
    \end{equation}

    We first establish the proof in the case $f\geq \lambda$.
    Let $g$ denote the pointwise maximal element of $\mathcal{G}(\lambda,\alpha)$ satisfying $g\leq f$.
    Such a function exists by \cref{p:relative-compact}.
    Suppose $x_0\in\R$ is arbitrary.
    For each $0\leq \delta \leq \alpha-g(x_0)$, let $h_\delta\in\mathcal{G}(\lambda,\alpha)$ denote the minimal function satisfying $h_\delta(x_0)=g(x_0)+\delta$.
    Equivalently, let $y_\delta$ be chosen so that
    \begin{equation*}
        \alpha - (\alpha - \lambda)\exp(y_\delta-x_0) = g(x_0)+\delta,
    \end{equation*}
    and define
    \begin{equation*}
        h_\delta(x)\coloneqq\begin{cases}
            \lambda &: x\leq y_\delta,\\
            \alpha - (\alpha - \lambda)\exp(y_\delta-x) &: y_\delta\leq x\leq x_0,\\
            \lambda - \bigl(\lambda - (g(x_0)+\delta)\bigr)\exp(x_0-x) &: x_0\leq x.
        \end{cases}
    \end{equation*}
    By minimality, $h_0\leq g$.
    Moreover, by taking $x_0$ sufficiently large, we may assume that $y_\delta\geq z$ for all $\delta>0$.

    If $g(x_0)=\alpha$, then the bound $f(x_0)\leq g(x_0)+e(x_0)$ is immediate from our assumption on $f$.
    Otherwise, if $g(x_0)<\alpha$, then for $\delta>0$, $g_\delta\coloneqq \max\{h_\delta, g\}\in\mathcal{G}(\lambda,\alpha)$ and $g_\delta(x_0)>g(x_0)$, so by maximality of $g$ there exists $y$ such that
    \begin{equation*}
        f(y)<h_\delta(y).
    \end{equation*}
    Since $h_\delta(x)\leq f(x)$ for all $x\leq y_\delta$, we must have $y\geq z$.
    If $y\leq x_0$ then
    \begin{align*}
        f(x_0) &\leq \alpha - (\alpha - f(y))\exp(y-x_0) + e(x_0)\\
         &< \alpha - (\alpha - h_\delta(y))\exp(y-x_0) + e(x_0)\\
         &\leq g(x_0)+\delta+e(x_0).
    \end{align*}
    Similarly if $y\geq x_0$ then
    \begin{equation*}
        \lambda - \bigl(\lambda - (g(x_0)+\delta)\bigr)\exp(x_0-y)> f(y)\geq \lambda - (\lambda - f(x_0) + e(x_0))\exp(x_0-y)
    \end{equation*}
    so $f(x_0)< g(x_0)+e(x_0)+\delta$.
    Since $\delta>0$ was arbitrary, it follows that if $g(x_0)<\alpha$ then
    \begin{equation*}
        f(x_0) \leq g(x_0) + e(x_0).
    \end{equation*}
    From this it follows that $\lim_{x\to\infty}\bigl(f(x)-g(x)\bigr)=0$.

    To establish the general case, one can apply the above strategy to the function $\max\{f,\lambda\}$, and the lower bound also follows by \cref{e:crude-bound}.
\end{proof}
We can now prove \cref{p:covering-class} and show that the covering class is well-defined.
In fact, we prove a slightly improved bound using the \emph{quasi-Assouad dimension} of $E$ \cite{zbl:1345.28019}, which is defined as follows:
\begin{align*}
    \dimqA E=\lim_{\theta\to 1^-}\inf\Bigl\{t\geq 0&:(\exists C>0)\,(\forall 0<r\leq R^{1/\theta}\leq R<1)\,(\forall x\in E)\\
                                                             & N_r\bigl(E\cap B(x,R)\bigr)\leq C\left(\frac{R}{r}\right)^t\Bigr\}.
\end{align*}
If $E\subset X$ where $X$ is Ahlfors--David $\alpha$-regular, then $\dimqA E\leq \alpha$.
In particular, if $E\subset\R^d$, then $0\leq\dimqA E\leq d$.
More detail on the quasi-Assouad dimension, and related Assouad-type dimensions, can be found in \cite{zbl:1467.28001}.
\begin{proposition}\label{p:covering-class-improved}
    Let $E\subset\R^d$ be non-empty and bounded with associated function $f\colon\R\to\R$ defined by
    \begin{equation*}
        f(x)=s_E(\exp(-\exp(x))).
    \end{equation*}
    Then there exists $g\in\mathcal{G}\bigl(0,\dimqA E\bigr)\subseteq \mathcal{G}(0,d)$ such that $f\eqinf g$.
\end{proposition}
\begin{proof}
    Write $\alpha=\dimqA E$ and set $\alpha_n = \alpha+1/n$.
    By the definition of quasi-Assouad dimension, for each $n\in\N$, get a constant $C_n'>0$ so that for all $x\in E$ and $0<r\leq R^{(n+1)/n}\leq R < 1$,
    \begin{equation*}
        N_r\bigl(E\cap B(x,R)\bigr)\leq C_n' \left(\frac{R}{r}\right)^{\alpha_n}.
    \end{equation*}
    Therefore there exists $C_n > 0$ such that
    \begin{equation*}
    N_r(E) \leq C_n N_R(E) \left( \frac{R}{r}\right)^{\alpha_n}.
    \end{equation*}
    Rearranging this bound, we obtain for $x_0\in\R$ and $x\geq \log\left(\frac{n+1}{n}\right)$,
    \begin{align*}
        f(x_0+x) &\leq \alpha_n - (\alpha_n- f(x_0))\exp(-x) + C_n\exp(-x_0-x)\\
                 &\leq \alpha - (\alpha - f(x_0))\exp(-x) + C_n\exp(-x_0-x) + \frac{1}{n}.
    \end{align*}
    Of course, we also have the bounds from the ambient Euclidean space as guaranteed by \cref{l:covering-reg}: for all $x_0\in\R$ and $x\geq 0$,
    \begin{equation*}
        f(x_0)\exp(-x)\leq f(x_0+x) \leq d - (d- f(x_0))\exp(-x) + A_d\exp(-x_0-x).
    \end{equation*}
    The upper bound further implies that for $0\leq x\leq \log\left(\frac{n+1}{n}\right)$,
    \begin{align*}
        f(x_0+x) \leq \alpha - (\alpha- f(x_0))\exp(-x) + A_d\exp(-x_0-x) + \frac{d-\alpha}{n+1}.
    \end{align*}
    Thus for $y\in\R$, set
    \begin{equation*}
        e(y)\coloneqq\inf_{n\in\N}\max\left\{A_d\exp(-y) + \frac{d-\alpha}{n+1}, C_n\exp(-y) + \frac{1}{n}\right\}.
    \end{equation*}
    Of course $e(y)\to 0$ as $y \to \infty$.
    Moreover, for all $x_0\in\R$ and $x\geq 0$,
    \begin{equation*}
        f(x_0+x) \leq \alpha - (\alpha- f(x_0))\exp(-x) + e(-x_0-x).
    \end{equation*}
    Thus the result follows by \cref{l:approx-bound}.
\end{proof}

\subsection{An alternative asymptotic formula}\label{ss:alt}
Finally, we can prove \cref{it:rate-form}.
In fact, it is a direct consequence of the following proposition combined with \cref{it:main}.
\begin{proposition}\label{p:char}
    Let $0\leq \lambda\leq\alpha$ be arbitrary and let $f\in\mathcal{G}(0,\alpha)$.
    Let
    \begin{equation*}
        g(x) = \sup_{\theta\in(0,1]} \left((1-\theta) \lambda + \theta f\bigl(x-\log(1/\theta)\bigr)\right).
    \end{equation*}
    Then $g$ is the pointwise minimal element of $\mathcal{G}(\lambda,\alpha)$ bounded below by $f$.

    Moreover, if $\diniu+g(x)\neq \lambda - g(x)$, then $g(x)=f(x)$ and $\diniu+g(x)=\diniu+f(x)$.
\end{proposition}
\begin{proof}
	Taking $\theta=1$, we see that $g \geq f$.
    Let $x_0\in\R$ be fixed, and let $\varepsilon>0$.
    By the definition of $g$, get $0<\theta_0\leq 1$ such that
    \begin{equation}\label{e:close-sup}
        (1-\theta_0) \lambda + \theta_0 f\bigl(x_0-\log(1/\theta_0)\bigr)\geq g(x_0)-\varepsilon.
    \end{equation}
    Now suppose $g_0\in\mathcal{G}(\lambda,\alpha)$ satisfies $g_0 \geq f$ pointwise.
    Then
    \begin{align*}
        g_0(x_0) &\geq \lambda - \left(\lambda - g_0\bigl(x_0-\log(1/\theta_0)\bigr)\theta_0 \right)\\
               &\geq (1-\theta_0)\lambda + \theta_0 f\bigl(x_0-\log(1/\theta_0)\bigr)\\
               &\geq g(x_0)-\varepsilon.
    \end{align*}
    Since $x_0\in\R$ and $\varepsilon>0$ were arbitrary, it follows that $g_0\geq g$.

    In the remainder of the proof, we show that $g\in\mathcal{G}(\lambda,\alpha)$.
    Let $x_0\in\R$, $x>0$, and $\varepsilon>0$ be arbitrary.
    Again by the definition of $g$, get $0<\theta_0\leq 1$ so that \cref{e:close-sup} holds and let $\theta=\exp(-x)$.
    Then
    \begin{align*}
        g(x_0+x)&\geq (1-\theta\theta_0)\lambda + \theta\theta_0 f\bigl(x_0+x-\log(1/(\theta\theta_0))\bigr)\\
            &= (1-\theta)\lambda + \theta\left((1-\theta_0)\lambda + \theta_0 f\bigl(x_0-\log(1/\theta_0)\bigr)\right)\\
            &\geq (1-\theta)\lambda + \theta(g(x_0)-\varepsilon)\\
            &\geq \lambda - (\lambda - g(x_0))\exp(-x) - \varepsilon.
    \end{align*}
    This gives the first inequality in \cref{l:scale-bound}, and also implies that $\diniu+g(x)\geq \lambda -g(x)$ for all $x\in\R$.

    To complete the proof, it suffices to establish the following fact: if $\diniu+ g(x) > \lambda - g(x)$, then $g(x)=f(x)$ and $\diniu+ g(x)=\diniu+ f(x)$.
    Suppose $x\in\R$ is such that $\diniu+ g(x)>\lambda -g(x)$.
    Equivalently, there exists $\kappa>\lambda$ and a sequence $(x_n)_{n=1}^\infty$ converging to $x$ from the right and $\kappa_n$ converging to $\kappa$ so that
    \begin{equation}\label{e:kappa-def}
        g(x_n) = \kappa_n - (\kappa_n - g(x))\exp(x-x_n).
    \end{equation}
    Now, suppose $y\leq x$ is arbitrary and let $\kappa>\kappa'>\lambda$.
    We may assume that $\kappa_n\geq\kappa'$ for all $n$.
    As proved above,
    \begin{equation*}
        g(x)\geq \lambda - (\lambda - g(y))\exp(y-x).
    \end{equation*}
    Combining the previous equations,
    \begin{align*}
        g(x_n) &\geq \kappa' - (\kappa' - (\lambda - (\lambda - g(y))\exp(y-x)))\exp(x-x_n)\\
               &\geq (1-\exp(y-x_n))\lambda + g(y)\exp(y-x_n) + \delta,
    \end{align*}
    where
    \begin{equation*}
        \delta\coloneqq (1-\exp(x-x_n))(\kappa'-\lambda)>0
    \end{equation*}
    is a constant which does not depend on $y$.
    In particular, by the definition of $g$, for each $n$ there exists $\theta_n$ such that $y_n = x_n-\log(1/\theta_n)\geq x$ and
    \begin{equation}\label{e:yn-equal}
        g(x_n) = (1-\theta_n) \lambda + \theta_n f(y_n).
    \end{equation}
    Since $x_n$ converges to $x$, $\theta_n$ converges to $1$ so by continuity of $f$, $g(x)=f(x)$.
    Finally, since $g(x)=f(x)$ and $f\leq g$, it remains to show that $\diniu+g(x)\leq\diniu+ f(x)$.
    To see this, again combining \cref{e:kappa-def} and \cref{e:yn-equal}, for each $n$,
    \begin{equation*}
        f(y_n) \geq \kappa_n - (\kappa_n - f(x))\exp(x-y_n).
    \end{equation*}
    But $\kappa_n$ converges to $\kappa$ and $y_n$ converges to $x$, so $\diniu+ f(x)\geq\diniu+g(x)$ as claimed.
\end{proof}
\begin{remark}
    Note that if $\diniu+ f(x)>\lambda-f(x)$, this \emph{does not} imply that $f(x)=g(x)$.
    For a visual depiction of this fact, see \cref{f:threshold}.
\end{remark}

\section{Examples and applications}\label{s:examples}
\subsection{Constructing countable discrete sets}
In this section, we will demonstrate the existence of countable discrete sets with various approximation properties.
To do so, we will use homogeneous Moran sets.
The construction of such sets is analogous to the usual $2^d$-corner Cantor set, except that the subdivision ratios need not be the same at each level.

Set $\mathcal{J}=\{0,1\}^d$.
We write $\mathcal{J}^*=\bigcup_{n=0}^\infty\mathcal{J}^n$, and we denote the word of length $0$ by $\varnothing$.
Suppose we have a sequence $\bm{r}=(r_n)_{n=1}^\infty$ with $0<r_n\leq 1/2$ for each $n\in\N$.
Then for all $n$ and $\bm{i}\in\mathcal{J}$, we define $S^n_{\bm{i}}\colon\R^d\to\R^d$ by
\begin{equation*}
    S^n_{\bm{i}}(x)\coloneqq r_n x+b^n_{\bm{i}}
\end{equation*}
where $b^n_{\bm{i}}\in\R^d$ is given by
\begin{equation*}
    (b^n_{\bm{i}})^{(j)} \coloneqq
    \begin{cases}
        0 &: \bm{i}^{(j)}=0\\*
        1-r_n &: \bm{i}^{(j)}=1
    \end{cases}.
\end{equation*}
We then set
\begin{equation*}
    M_n\coloneqq \bigcup_{(\bm{i}_1,\ldots,\bm{i}_n)\in\mathcal{J}^n}S^1_{\bm{i}_1}\circ\cdots\circ S^n_{\bm{i}_n}([0,1]^d)\qquad\text{and}\qquad M=M(\bm{r})\coloneqq\bigcap_{n=1}^\infty M_n.
\end{equation*}
We call $M$ a \emph{homogeneous Moran set}.
Note that $M_n$ consists of $2^{dn}$ hypercubes, each with side-length $\rho_n\coloneqq r_1\cdots r_n$.

We first demonstrate the existence of homogeneous Moran sets with arbitrary covering class.
\begin{lemma}\label{l:existmoran}
    Let $d\in\N$ and $g\in\mathcal{G}(0,d)$ be arbitrary.
    Then there exists a homogeneous Moran set with covering class $g$.
\end{lemma}
\begin{proof}
    First assume that $g \not\eqinf 0$, and note that $\limsup_{x \to \infty} g(x) > 0$.
    Therefore for all $y\in\R$ there exists a minimal $\psi(y)>y$ so that
    \begin{equation*}
        g(y)\exp(y-\psi(y))=g(\psi(y))-d\log(2)\cdot\exp(-\psi(y)).
    \end{equation*}
    Now set $x_1=0$ and, inductively, set $x_{k+1}=\psi(x_k)$ for each $k\in\N$.
    Let $\rho_k=\exp(-\exp(x_k))$ denote the corresponding scales (note that $\rho_1 = r_1$), and set $r_k\coloneqq\rho_k/\rho_{k-1}$ for $k\geq 2$.
    Thus for $0<\delta\leq r_1$, if $k$ is such that $\rho_k<\delta\leq\rho_{k-1}$, we set
    \begin{equation*}
        \overline{s}(\delta)=\frac{kd\log 2}{\log(1/\delta)}.
    \end{equation*}
    The exact same calculation as in the proof of \cite[Lemma~3.4]{zbl:1509.28005} gives that for all $k\geq 2$ we have $r_k\in(0,1/2]$, $\overline{s}(\rho_k)=g(x_k)$, and
    \begin{equation}\label{e:bound}
        g(x)-d\log(2)\exp(-x)\leq\overline{s}(\exp(-\exp(x)))\leq g(x)
    \end{equation}
    for all $x \geq x_2$.
    In particular, the resulting homogeneous Moran set has covering class $g$.

    Finally, if $g \eqinf 0$, then it is straightforward to check directly that the homogeneous Moran set given by sequence $r_n = 2^{-2^n}$ has covering class $g$.
\end{proof}

Now by discretising a homogeneous Moran set, we can obtain a countable discrete set with arbitrary covering class.
\begin{lemma}\label{l:countable}
	Let $d\in\N$ and $g\in\mathcal{G}(0,d)$ be arbitrary.
	Then there exists a countable discrete set $F \subset (0,1)^d$ which accumulates only at $0$ and has covering class $g$.
\end{lemma}
\begin{proof}
	By \cref{l:existmoran}, get a homogeneous Moran set $M$ with covering class $g$.
	For all $n \in \N$ let $F_n$ be a finite subset of $M \cap (0,1/n)^d$ whose Hausdorff distance from $M \cap [0,1/n]^d$ is at most $2^{-n}$.
	Then
    \begin{equation*}
        F \coloneqq \bigcup_{n=1}^{\infty} F_n
    \end{equation*}
    is clearly discrete and accumulates only at $0$.
	If $2^{-n} \leq r < 2^{-(n-1)}$, then
    \begin{equation*}
        N_r(F) \geq N_r(F_n) \approx N_r(M \cap [0,1/n]^d) \approx n^{-d} N_r(M),
    \end{equation*}
	with implicit constants independent of $n$.
	In particular, as $x \to \infty$,
    \begin{equation*}
        | s_F(\exp(-\exp(x)))  - s_M(\exp(-\exp(x)))| \lesssim xe^{-x} \xrightarrow[]{} 0.
    \end{equation*}
    Therefore since $M$ has covering class $g$, $F$ also does.
\end{proof}
Finally, we use \cref{l:countable} to obtain an infinitely generated self-similar IFS with specified Hausdorff dimension and fixed points having arbitrary covering class.
\begin{lemma}\label{l:fixed-point-class}
	Let $d\in\N$.
    Let $g\in\mathcal{G}(0,d)$ and $0<h<d$ be arbitrary.
    Then there exists a countable self-similar IFS $\{S_i\}_{i\in\mathcal{I}}$ with fixed points $F$ and attractor $\limitset$ such that $F$ has covering class $g$ and $\dimH\limitset=h$.
\end{lemma}
\begin{proof}
    Since $h<d$, we may fix $n$ large enough that
    \begin{equation}\label{e:number-of-hypercubes}
        \frac{\log (2^{dn} - 1)}{n\log 2} > h.
    \end{equation}
    By \cref{l:countable}, let $F_0\subset(0,1)^d$ be a countable discrete set which accumulates only at $0$ and has covering class $g$, and let $F_1\coloneqq F_0\cap (0,2^{-n})^d$.
    To each $p \in F_1$ we will now choose a similarity map $S_p$ which fixes $p$ and has some contraction ratio $c_p \in (0,1)$.
    We choose the contraction ratios to be small enough that $S_p([0,1]^d) \subset  (0,2^{-n})^d$ for all $p \in F_1$, and $S_p([0,1]^d) \cap S_q([0,1]^d) = \varnothing$ whenever $p,q \in F_1$ are distinct, and moreover $\sum_{p \in F_1}c_p^h < 1$.
    By \cref{e:number-of-hypercubes} there exists $c \in (0,2^{-n})$ such that
    \begin{equation}\label{e:hausdorff-choice}
        (2^{dn} - 1) c^h = 1-\sum_{p \in F_1}c_p^h.
    \end{equation}
    Fix similarity maps $T_1,\dotsc,T_{2^{nd}-1}$, each with contraction ratio $c$, such that whenever $1 \leq i < j \leq 2^{nd}-1$, we have that $T_i([0,1]^d)$ and $T_j([0,1]^d)$ are pairwise-disjoint subsets of $(0,1)^d \setminus (0,2^{-n}]^d$.

    Now consider the countable self-similar IFS
    \begin{equation*}
    \{S_p : p \in F_1\} \cup \{T_i\}_{1 \leq i \leq 2^{nd}-1}.
    \end{equation*}
    Since $F_0$ accumulates only at $0$, the symmetric difference of $F_0$ and the set $F$ of fixed points of this \cifs{} is finite.
    Therefore since $F_0$ has covering class $g$, the same is true for $F$.
    Moreover, combining \cref{e:hausdorff-formula} due to Mauldin \& Urbański with \cref{e:hausdorff-choice}, the Hausdorff dimension of the limit set equals $h$.
\end{proof}

Finally, using the construction established above, we prove \cref{it:fixed}.
For convenience of notation, we make one more definition.
\begin{definition}
    Given a sequence of functions $(f_k)_{k=1}^\infty$ each defined on some interval $[0,a_k]$, the \defn{concatenation} of $(f_k)_{k=1}^\infty$ is the function $f\colon(-\infty,\sum_{k=1}^\infty a_k)\to\R$ given as follows: for each $x>0$ with $\sum_{j=0}^{k-1} a_j<x\leq\sum_{j=0}^{k}a_j$ where $a_0=0$ we define
    \begin{equation*}
        f(x)=f_k\left(x-\sum_{j=0}^{k-1} a_j\right),
    \end{equation*}
    and for $x\leq 0$ we define $f(x)=f_1(0)$.
\end{definition}
\begin{figure}[t]
    \centering
    \begin{tikzpicture}[>=stealth,xscale=6,yscale=9]
    \begin{scope}[dotted]
        \draw (-0.1,0.26)node[left]{$h$} -- (2.2,0.26);
        \draw (-0.1,0.3)node[left]{$s$} -- (2.2,0.3);
        \draw (-0.1,0.369)node[left]{$\beta$} -- (2.2,0.369);
        \draw (-0.1,0.7)node[left]{$t$} -- (2.2,0.7);
    \end{scope}

    \begin{scope}[dotted]
        \draw (0.847298,0.26-0.03) -- (0.847298,0.7+0.03);
        \draw (1.29165,0.26-0.03) -- (1.29165,0.7+0.03);
        \draw (1.39543,0.20-0.03) -- (1.39543,0.7+0.03);
    \end{scope}

    \draw[<->] (0, 0.26) -- node[fill=white]{$a_{1,n}$} (0.847298, 0.26);
    \draw[<->] (0.847298, 0.26) -- node[fill=white]{$a_{2,n}$} (1.29165, 0.26);
    \draw[<->] (1.29165, 0.26) -- node[fill=white]{$a_{3,n}$} (2.13895, 0.26);

    \draw[<->] (0, 0.20) -- node[fill=white]{$b_{1,n}$} (1.39543, 0.20);
    \draw[<->] (1.39543, 0.20) -- node[fill=white]{$b_{2,n}$} (2.13895, 0.20);

    \draw[thick, dashed] plot file {figures/sharpness/original.txt};
    \draw[thick] plot file {figures/sharpness/modified.txt};
\end{tikzpicture}
    \caption{A plot of the covering class of $F$ (dashed) and the covering class of $\limitset$ (solid) corresponding to the concatenation of $(g_{1,n},g_{2,n})$ and $(f_{1,n},f_{2,n},f_{3,n})$ respectively.
        In this plot, we assume that $h<s<\beta<t$ to remove the dependence on $n$.
    }
    \label{f:sharpness}
\end{figure}
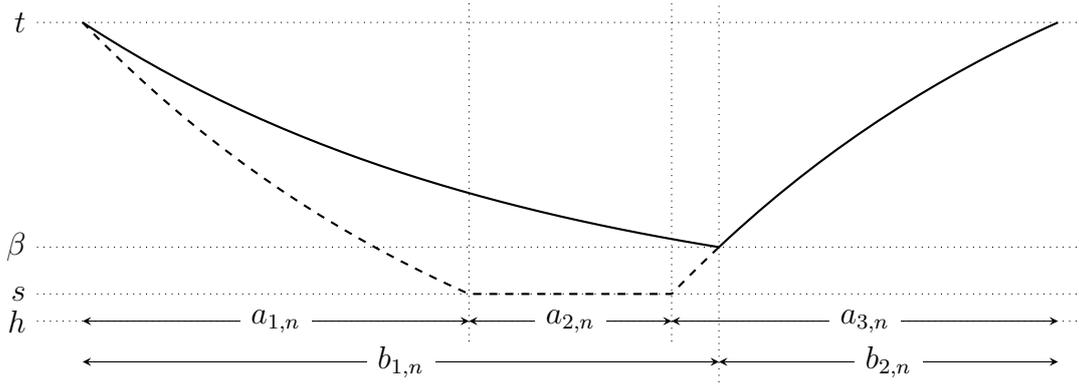
\begin{proofref}{it:fixed}
    The bounds on $\dimlB\limitset$ have already been proved in \cref{c:fixed-bound}; what remains is to verify sharpness.

    We first observe that the result is straightforward to prove if $\mathcal{D}(h,s,t,d)$ is a singleton.
    If this is the case, let $g\in\mathcal{G}(0,d)$ be such that $\liminf_{x\to \infty}g(x)=s$ and $\limsup_{x\to\infty}g(x)=t$.
    Applying \cref{l:fixed-point-class}, get a self-similar IFS $\{S_i\}_{i\in\mathcal{I}}$ with fixed points $F$ and attractor $\limitset$ such that $F$ has covering class $g$.
    Then by \cref{c:fixed-bound}, since $\mathcal{D}(h,s,t,d)$ is a singleton, it must be that $\dimlB\limitset$ is the expected value.

    Otherwise, $\mathcal{D}(h,s,t,d)$ is not a singleton, so $0<h<t$ and $0<s<t$.
    Let $\beta\in\mathcal{D}(h,s,t,d)$, or equivalently
    \begin{equation}\label{e:beta-bound}
        \max\{s,h\}\leq \beta \leq h + \frac{(t-h) (d - h) s}{d\cdot t-h\cdot s}.
    \end{equation}
    Given $\delta>0$, set
    \begin{equation*}
        \beta_n\coloneqq \max\left\{\beta, h+\frac{\delta}{n}\right\}
        \qquad\text{and}\qquad
        t_n\coloneqq \min\left\{t, d - \frac{\delta}{n}\right\}.
    \end{equation*}
    Note that since $t>\max\{h,s\}$, it follows that $h\leq \beta < t$, so by taking $\delta>0$ sufficiently small we may assume that
    \begin{equation}\label{e:delta-small}
        h<\beta_n<t_n<d\qquad\text{and}\qquad \beta_n\leq h + \frac{(t-h) (d - h) s}{d\cdot t-h\cdot s}.
    \end{equation}
    (These choices are only actually necessary in the case that $\beta=h$ or $t=d$; otherwise, we could just take $\beta_n=\beta$ and $t_n=t$ for all $n\in\N$.)

    We now choose some constants.
    Let $a_{1,n}>0$ be chosen so that $t_n\exp(-a_{1,n})=s$ and let $a_{3,n}>0$ be chosen so that $d-(d-s)\exp(-a_{3,n})=t_{n+1}$.
    Then let $a_{3,n}' \in (0,a_{3,n}]$ be chosen so that $d-(d-s)\exp(-a_{3,n}')=\beta_n$.
    Finally, let $a_{2,n}'\in \R$ be chosen so that
    \begin{equation*}
        h-(h-t_n)\exp(-a_{1,n}-a_{2,n}'-a_{3,n}')=\beta_n,
    \end{equation*}
    and let $a_{2,n} \coloneqq \max\{a_{2,n}',0\}$.
    Note that $a_{2,n}'\geq 0$ if and only if the second equation in \cref{e:delta-small} holds with $t_n$ in place of $t$, so since $t_n \to t$ we have $\liminf_{n \to \infty} a_{2,n}' \geq 0$.
    Note also that $0 < a_{1,1} \leq a_{1,2} \leq a_{1,3} \leq \dotsb$ and $0 < a_{3,1} \leq a_{3,2} \leq a_{3,3} \leq \dotsb$.

    Now, we define functions
    \begin{itemize}[nl]
        \item $f_{1,n}(x) = t_n\exp(-x)$ for $x\in[0,a_{1,n}]$;
        \item $f_{2,n}(x)=s$ for $x\in[0,a_{2,n}]$; and
        \item $f_{3,n}(x)=d-(d-s)\exp(-x)$ for $x\in[0,a_{3,n}]$.
    \end{itemize}
    Let $f$ be the concatenation of the sequence of functions
    \begin{equation*}
        (f_{1,1},f_{2,1},f_{3,1},f_{1,2},f_{2,2},f_{3,2},f_{1,3},\ldots).
    \end{equation*}
    Of course, $f\in\mathcal{G}(0,d)$ and moreover by construction $\liminf_{x\to \infty}f(x)=s$ and $\limsup_{x\to\infty}f(x)=\lim_{n\to\infty}t_n = t$.
    Moreover, the minimal function $g\in\mathcal{G}(h,d)$ satisfying $f\leq g$ is similarly the concatenation of the sequence of functions
    \begin{equation*}
        (g_{1,1},g_{2,1},g_{1,2},g_{2,2},g_{1,3},\ldots)
    \end{equation*}
    where, setting $b_{1,n} = a_{1,n}+a_{2,n}+a_{3,n}'$ and $b_{2,n}=a_{3,n}-a_{3,n}'$,
    \begin{itemize}[nl]
        \item $g_{1,n}(x)=h-(h-t_n)\exp(-x)$ for $x\in[0,b_{1,n}]$; and
        \item $g_2(x) = (d-(d-\beta_n)\exp(-x)$ for $x\in [0,b_{2,n}]$.
    \end{itemize}
    A depiction of the functions $f$ and $g$ can be found in \cref{f:sharpness}.

    To conclude the proof, by \cref{l:fixed-point-class}, get a self-similar IFS $\{S_i\}_{i\in\mathcal{I}}$ with fixed points $F$ and attractor $\limitset$ such that $F$ has covering class $g$ and $\dimH\limitset = h$.
    Then recalling \cref{it:rate-form},
    \begin{equation*}
        \dimlB\limitset = \liminf_{x\to\infty}g(x)=\lim_{n\to\infty}\beta_n=\beta
    \end{equation*}
    as required.
\end{proofref}

\subsection{Continued fraction expansions with restricted entries}\label{ss:cont-frac}
In this section we prove \cref{it:gauss-invariant}.
For a non-empty, proper subset $I \subset \N$, define
\begin{equation*}
    \limitset_I \coloneqq \left\{ z \in (0,1) \setminus \mathbb{Q} : z = \frac{1}{b_1 + \frac{1}{b_2 + \frac{1}{\ddots}}}, b_n \in I \mbox{ for all } n \in \N \right\}.
\end{equation*}
It is well-known (see, for instance, \cite[p.~4997]{zbl:0940.28009}) that $\limitset_I$ is the limit set of the \cifs{} given by the inverse branches of the Gauss map corresponding to the elements of $I$.
Indeed, this is one of the motivations for working with countable IFSs given by conformal maps rather than just similarity maps.
\begin{lemma}\label{l:ctdfraccifs}
    Working in $\mathbb{R}$, letting $X \coloneqq [0,1]$,
    \begin{itemize}[nl]
        \item\label{i:ctdfracnot1} If $1 \notin I$ then $\{ S_b(x) \coloneqq 1/(b+x) : b \in I \}$ is a \cifs{} with limit set $\limitset_I$.
        \item\label{i:ctdfrac1} If $1 \in I$ then $\{ S_b(x) \coloneqq 1/(b+x) : b \in I, b \neq 1 \} \cup \left\{ S_{1b}(x) \coloneqq \frac{1}{b+\frac{1}{1+x}} : b \in I \right\}$ is a \cifs{} with limit set $\limitset_I$.
    \end{itemize}
\end{lemma}
We can finally prove our headline result.
In fact, since continued fraction sets $\Lambda_I$ are Borel and invariant for the Gauss map, the following stronger result immediately implies \cref{it:gauss-invariant}.
\begin{theorem}
    For all $0<a<1$ there exists an infinite proper subset $I \subset \N$ such that the box dimension of the continued fraction set $\Lambda_I$ does not exist, and $\Lambda_I\subset(0,a)$.
\end{theorem}
\begin{proof}
	We define a sequence $(a_n)_{n \geq 0}$ inductively by setting $a_0 = 2$ and $a_n = (2a_{n-1})^n$ for $n \geq 1$.
    Then let
    \begin{equation*}
        I_0 = \left\{b^2 : b \in \N \cap \bigcup_{n=0}^{\infty} [a_n,2a_n] \right\}.
    \end{equation*}
    Now, using notation from \cref{l:ctdfraccifs}, we have $|S_b'(x)| \approx b^{-2}$ uniformly for $b \geq 2$ and $x \in [0,1]$.
    Therefore
    \begin{equation*}
        \pressure(t)<\infty\qquad\text{for all}\qquad t>1/4.
    \end{equation*}
    In particular, by \cite[Theorem~3.23]{zbl:0852.28005} there exists $N \in \N$ large enough that if
    \begin{equation*}
        I \coloneqq I_0 \cap [N,\infty)
    \end{equation*}
    then $\dimH \limitset_I < 1/3$.
    We may increase $N$ further if necessary so that $\limitset_I\subset(0,a)$.

    Now, note that the orbit set $Q_I \coloneqq \mathcal{O}(0,1) = \{1/b : b \in I_0\}$ is a discrete approximation of the \cifs{} $\{ S_b : b \in I_0 \}$.
    Let $\fixedpts_I$ be the set of fixed points and let $\Lambda_I$ be the limit set.
    A mean value theorem argument gives that $(b+1)^{-2} - b^{-2} \approx b^{-3}$.
    Therefore using \cref{l:fixed-choice},
    \begin{equation*}
        N_{(2a_n)^{-3}}(\fixedpts_I) \approx N_{(2a_n)^{-3}}(Q_I) \geq N_{(2a_n)^{-3}}(\{b^{-2} : a_n \leq b \leq 2a_n\}) \gtrsim a_n.
    \end{equation*}
    In particular, $\dimuB \fixedpts_i \geq 1/3$.
    On the other hand, for $n \geq 1$,
    \begin{equation*}
        N_{a_n^{-2}}(\fixedpts_I) \approx N_{a_n^{-2}}(Q_I) \leq 2 + a_{n-1} \lesssim a_n^{1/n}.
    \end{equation*}
    This implies that $\dimlB \fixedpts_I \leq 1/n$ for all $n$, so $\dimlB \fixedpts_I = 0$.
    In particular,
    \begin{equation*}
        \max \{\dimH \limitset_I, \dimlB \fixedpts_I\} = \dimH \limitset_I < 1/3 \leq \dimuB \fixedpts_I.
    \end{equation*}
    Thus by \cref{it:box-exist}, the box dimension of $\limitset_I$ does not exist.
\end{proof}

\begin{acknowledgements}
	We thank Jonathan Fraser, Mike Todd, and Mariusz Urbański for discussions related to the contents of this paper.
	We thank Simon Baker for comments on a draft version of this manuscript and Thomas Jordan for suggesting some relevant references.
    AB was supported by an EPSRC New Investigators Award (EP/W003880/1).
    AR was supported by EPSRC Grant EP/V520123/1 and the Natural Sciences and Engineering Research Council of Canada.
\end{acknowledgements}
\end{document}